\documentclass[12pt, a4paper]{amsart}

\usepackage[hmargin=30mm, vmargin=25mm, includefoot, twoside]{geometry}
\usepackage[
  bookmarksopen=true,
  hypertexnames=false]{hyperref}

\usepackage{amsfonts,amssymb,verbatim}

\usepackage{latexsym}
\usepackage{mathrsfs}
\usepackage{stmaryrd}
\usepackage{xspace}
\usepackage{enumerate, paralist}
\usepackage{graphicx}
\usepackage[all]{xy}
\usepackage{extarrows}
\usepackage{tabu}
\usepackage{accents}
\usepackage{tensor}

\usepackage[usenames,dvipsnames]{color}

\usepackage{txfonts, pxfonts}

\usepackage{amsmath}

\usepackage[tbtags]{mathtools}

\usepackage{amsthm, thmtools}

\newtheorem{thm}{Theorem}[section]

 \newtheorem{cor}[thm]{Corollary}
 \newtheorem{lem}[thm]{Lemma}
 \newtheorem{prop}[thm]{Proposition}

\newtheorem{introthm}{Theorem}

\newtheorem{introcor}[introthm]{Corollary}
\numberwithin{equation}{section}

 \theoremstyle{definition}
  \newtheorem{defn}[thm]{Definition}

 \theoremstyle{remark}
 \newtheorem{rem}[thm]{Remark}

\newtheorem*{claim*}{Claim}

\def\NN{\mathbb{N}}
\def\CC{\mathbb{C}}

\def\ZZ{\mathbb{Z}}

\def\B{\mathfrak{B}}

\def\KK{\mathfrak{K}}

\def\H{\mathcal{H}}
\def\K{\mathcal{K}}
\def\N{\mathcal{N}}

\def\G{\mathcal{G}}

\def\MM{\mathcal{M}}
\def\S{\mathcal{S}}

\def\Id{\mathrm{Id}}
\def\I{\mathrm{I}}
\def\Tr{\mathrm{Tr}}

\def\ppg{\mathrm{prop}}

\def\dim{\mathrm{dim}}

\def\rank{\mathrm{rank}}
\def\supp{\mathrm{supp}}
\def\dist{\mathrm{dist}}
\def\diam{\mathrm{diam}}

\def\girth{\mathrm{girth}}

\def\deg{\mathrm{deg}}
\def\mdeg{\mathrm{m\text{-}deg}}

\begin{document}

\title{Quasi-locality and weak containment for Roe algebras}
\author{Benyin Fu}
\author{Kang Li}
\author{Jiawen Zhang}

\address[Benyin Fu]{College of Statistics and Mathematics, Shanghai Lixin University of Accounting and Finance, Shanghai, 201209, China.}
\email{fuby@lixin.edu.cn}

\address[Kang Li]{Centre for Mathematical Sciences, Lund University, Box 118, SE-221 00 Lund, Sweden.}
\email{kang.li@math.lth.se}

\address[Jiawen Zhang]{School of Mathematical Sciences, Fudan University, 220 Handan Road, Shanghai, 200433, China.}
\email{jiawenzhang@fudan.edu.cn}

\thanks{Benyin Fu was partly supported by Natural Science Foundation of Shandong Province (No. ZR2024MA048). Jiawen Zhang was partly supported by the National Key R{\&}D Program of China 2022YFA1007000 and NSFC 12422107.}


\begin{abstract}
We establish two characterisations of Property A for discrete metric spaces of bounded geometry in terms of the Roe algebra. More precisely, we show that such a space has Property A if and only if its Roe algebra coincides with its quasi-local algebra, also equivalent to the canonical quotient map from the maximal Roe algebra to the reduced one being injective. We also prove, without assuming Property A, that the strongly quasi-local algebra always coincides with the Roe algebra. As applications, we show that a discrete group is exact whenever its full group $C^*$-algebra has the local lifting property. We also show that the maximal injective crossed-product functor agrees with the reduced one precisely for exact discrete groups.
\end{abstract}

\date{}
\maketitle

\section{Introduction}
\label{sec:intro}

Roe algebras were introduced by John Roe in his work of index theory on open
manifolds \cite{Roe88}, where he showed that their $K$-theory provides a natural receptacle for
higher indices of elliptic operators. This insight has led to significant applications in geometry and topology \cite{Roe93, Roe96}. Since then, Roe algebras have become fundamental
objects in coarse geometry and higher index theory \cite{WillettYu2020}, connecting coarse 
properties of metric spaces with operator-algebraic invariants. 
However, it is often difficult to determine whether a given operator belongs to a Roe algebra, since its definition involves norm approximation by finite-propagation operators and there is no canonical procedure to construct such approximations. This leads to the question of whether other natural descriptions or
completions of the finite-propagation algebra give the same $C^*$-algebra.

In this paper, we investigate two such comparisons. The first replaces
finite propagation by quasi-locality, a notion introduced by
Roe \cite{Roe88} that controls the norm of compressions
between distant subsets. 
The second concerns the maximal Roe algebras, introduced in \cite{GWY08} to study the strong Novikov conjecture for residually finite groups.
The main geometric condition governing both comparisons
is Yu's Property~A, introduced in his study of the
coarse Baum-Connes conjecture \cite{Yu00}.
Property~A can be regarded as a coarse analogue of
amenability. For countable discrete groups, Property~A
is equivalent to exactness
(\cite{GK02, Ozawa2000}, see also \cite[Chapter~5]{BO}). Along the way, we also determine the
precise relationship between Roe algebras and the strongly quasi-local
algebras introduced in \cite{BCZ23, CGZ24}.

Throughout the paper, let $(X,d)$ be a discrete metric space of
bounded geometry (\emph{i.e.}, $\sup_{x\in X}|B(x,R)|<\infty$ for each $R>0$, where
$B(x,R)=\{y\in X:d(x,y)\le R\}$). Let $\H_0$ be a separable Hilbert
space and put $\H_X:=\ell^2(X;\H_0)\cong\ell^2(X)\otimes\H_0$.
For a subset $A\subseteq X$, write $\chi_A$ for both its characteristic
function and the corresponding multiplication projection. If
$T\in\B(\H_X)$ and $x,y\in X$, let
$T_{x,y}:=\chi_{\{x\}}T\chi_{\{y\}}\in\B(\H_0)$.
An operator $T$ has \emph{finite propagation} if there exists $R>0$
such that $T_{x,y}=0$ whenever $d(x,y)>R$. It is \emph{locally compact}
if $\chi_{\{x\}}T$ and $T\chi_{\{x\}}$ are compact for all $x\in X$.
The \emph{Roe algebra} $C^*(X;\H_0)$ is the norm closure in
$\B(\H_X)$ of the locally compact finite-propagation operators.
When $\H_0=\CC$, this is the usual uniform Roe algebra, also denoted by $C^*_u(X)$. When $\H_0$ is infinite-dimensional, this is the usual Roe algebra, also denoted by $C^*(X)$.

\subsection{Quasi-locality and strong quasi-locality}

Recall that $T\in\B(\H_X)$ is \emph{quasi-local} if for every
$\varepsilon>0$ there is $R>0$ such that $\|\chi_A T\chi_B\|<\varepsilon$ whenever $A,B\subseteq X$ with $d(A,B)>R$. 
The locally compact quasi-local operators form the
\emph{quasi-local algebra} $C_q^*(X;\H_0)$. We call $C^*_{uq}(X):=C_q^*(X;\CC)$ the \emph{uniform quasi-local algebra}. 
Clearly, every finite-propagation
operator is quasi-local, and hence 
$C^*(X;\H_0)\subseteq C_q^*(X;\H_0)$.
It is natural to ask when the two algebras coincide. \v{S}pakula and the third-named author \cite{SZ20} proved that $C^*(X;\H_0)=C_q^*(X;\H_0)$ when $X$ has Yu's Property~A, while Ozawa \cite{Oza} used expander graphs to show that they differ in general.

We give the following complete characterisation, answering \cite[Question 6.7]{SZ20} and confirming one of the
expectations in \cite{Oza}.

\begin{introthm}\label{introthm:uniform-converse}
Let $X$ be a discrete metric space of bounded geometry and $\H_0$
be a nonzero separable Hilbert space. Then $X$ has Property~A
\emph{if and only if}
$C^*(X;\H_0)=C_q^*(X;\H_0)$.
\end{introthm}

Combining Theorem~A and \cite[Proposition~6.1]{LNSZ21}, we obtain an affirmative answer to \cite[Question~7.1]{LNSZ21} as follows (see the definition therein).

\begin{introcor}
Let $X$ be a discrete metric space of bounded geometry.
Then $X$ has Property A \emph{if and only if} $\ell^\infty(X)$ is a Cartan subalgebra of $C_{uq}^*(X)$.
\end{introcor}

As a stronger version of the converse direction in Theorem~\ref{introthm:uniform-converse}, we show in Proposition~\ref{prop:ghost} that if $X$ does not have Property~A, then there is a ghost operator in $C^*_{uq}(X) \setminus C^*_u(X)$. An operator $T \in \B(\ell^2(X))$ is called a \emph{ghost} (a notion due to G. Yu) if for any $\varepsilon > 0$ there is a finite subset $F \subseteq X$ such that for any $(x,y) \notin F \times F$, we have $|T_{x,y}| < \varepsilon$. Together with \cite[Theorem 3.3]{SZ20}, this gives the following.

\begin{introcor}
Let $X$ be a discrete metric space of bounded geometry.
Then $X$ has Property A \emph{if and only if} any ghost in $C^*_{uq}(X)$ belongs to $C^*_u(X)$.
\end{introcor}

Theorem~\ref{introthm:uniform-converse} also has a groupoid
interpretation. Recall that Skandalis, Tu and Yu \cite{STY02} introduced the notion of coarse groupoid $G(X)$
for a discrete metric space $X$ of bounded geometry, whose topological amenability is
equivalent to Property~A of $X$. The notion of the quasi-local algebra for \'{e}tale groupoids was recently introduced in \cite{JZZ23} and it was proved that for locally compact, $\sigma$-compact amenable
\'{e}tale groupoids, the equivariant quasi-local algebra coincides
with the reduced groupoid $C^*$-algebra. Examples in
\cite{LZZ26} show that this equality can fail in general.
However, Theorem~\ref{introthm:uniform-converse} shows that for coarse
groupoids, this equality holds \emph{if and only if} the groupoid is
topologically amenable, since the groupoid construction recovers the uniform Roe and
uniform quasi-local algebras (see \cite{JZZ23}).

On the other hand, quasi-local algebras also arise naturally in higher index theory. Indeed, Engel showed that the higher indices of elliptic uniform pseudodifferential operators on open manifolds naturally take values in the $K$-theory of quasi-local algebras \cite{Eng18}. 
This motivates the finer question of whether the inclusion of
the Roe algebra into the quasi-local algebra induces an
isomorphism on $K$-theory. Recent work \cite{LZZ26} shows that
this map need not be surjective on $K_0$-groups for certain uniform
Roe algebras associated with large-girth expanders, while
the analogous comparison for the ordinary Roe and quasi-local
algebras remains open in general.

One potential approach to this $K$-theoretic problem is through strong quasi-locality, introduced in \cite{BCZ23, CGZ24}. It was proved in \cite{BCZ23} that the inclusion of the Roe algebra into the strongly quasi-local algebra induces an isomorphism on $K$-theory whenever $X$ admits a coarse embedding into Hilbert space. The definition is based on the commutator characterisation of ordinary quasi-locality by scalar-valued functions \cite[Theorem~2.8]{ST19}, by replacing scalar-valued test functions with compact-operator-valued ones.

More precisely, let us fix an infinite-dimensional
separable Hilbert space $\K$, and let $\KK(\K)_1$ be the closed unit
ball of the compact operators on $\K$. For $\delta,R>0$, a map
$g:X\to\KK(\K)_1$ has \emph{$(\delta,R)$-variation} if
$\|g(x)-g(y)\|<\delta$ whenever $d(x,y)<R$. Define
$\Lambda(g)\in\B(\H_X\otimes\K)$ by
\[
 \Lambda(g):=\bigoplus_{x\in X}(\Id_{\H_0}\otimes g(x))
\]
with respect to the decomposition
$\H_X\otimes\K\cong\bigoplus_{x\in X}(\H_0\otimes\K)$.
An operator $T\in\B(\H_X)$ is \emph{strongly quasi-local} if, for
every $\varepsilon>0$, there exist $\delta,R>0$ such that
\begin{equation}\label{EQ:commut strong quasi-local}
 \big\|[T\otimes\Id_{\K},\Lambda(g)]\big\|<\varepsilon
\end{equation}
for every map $g:X\to\KK(\K)_1$ with $(\delta,R)$-variation, where
$[A,B]:=AB-BA$. The \emph{strongly quasi-local algebra}
$C^*_{sq}(X;\H_0)$ consists of the locally compact strongly
quasi-local operators. By \cite[Proposition~3.6]{BCZ23}, $C^*(X;\H_0)\subseteq C^*_{sq}(X;\H_0)\subseteq C_q^*(X;\H_0)$.

Our next result identifies this intermediate algebra
without assuming Property~A or coarse embeddability.

\begin{introthm}\label{introthm:sq=Roe}
Let $X$ be a discrete metric space of bounded geometry and $\H_0$
be a separable Hilbert space. An operator
$T\in\B(\ell^2(X;\H_0))$ is strongly quasi-local \emph{if and only if}
for any $\varepsilon>0$, there exists a finite-propagation operator
$A\in\B(\ell^2(X;\H_0))$ such that $\|T-A\|\le\varepsilon$.
Furthermore, $C^*(X;\H_0)=C^*_{sq}(X;\H_0)$.
\end{introthm}

Thus, strong quasi-locality always characterises norm approximation
by finite-propagation operators, while ordinary quasi-locality does
so precisely under Property~A. 
Moreover, the inclusion $C^*(X;\H_0)\hookrightarrow C^*_{sq}(X;\H_0)$
induces an isomorphism on $K$-theory for every $X$, not only for
coarsely embeddable spaces.

\subsection{Maximal Roe algebras and weak containment}

Our second comparison concerns a different completion of the same
algebra. Let $\CC[X;\H_0]$ be the algebraic Roe algebra, consisting of
locally compact finite-propagation operators on $\H_X$. For $T\in\CC[X;\H_0]$, define
\[
\|T\|_{\max}:=\sup_{\pi} \|\pi(T)\|,
\]
where the supremum is taken over all $*$-representations
$\pi$ of $\CC[X;\H_0]$ on Hilbert spaces.
By \cite[Lemma~1.10]{SW}, $\|T\|_{\max}$ is always finite. The maximal Roe algebra $C^*_{\max}(X;\H_0)$ is the
completion of $\CC[X;\H_0]$ in this norm, and we denote the canonical quotient by 
$\lambda_{X,\H_0}:C^*_{\max}(X;\H_0)\to C^*(X;\H_0)$. When $\H_0=\CC$, we write $\CC_u[X]:=\CC[X;\CC]$, $C^*_{u,\max}(X):=C^*_{\max}(X;\CC)$ and $\lambda_X:C^*_{u,\max}(X)\to C_u^*(X)$.

\v{S}pakula--Willett \cite[Proposition~1.3]{SW} proved that these
quotients are isomorphisms when $X$ has Property~A. The converse was
left open in \cite[Section~1]{SW}.
In terms of coarse groupoids, this asks whether the
\emph{weak containment property}, namely the equality
of the full and reduced groupoid $C^*$-algebras,
implies topological amenability.
Willett \cite{Willett} showed that this fails
for general \'{e}tale groupoids, but the question remained open
for coarse groupoids (see \cite[Remark~3.5]{Willett}).


We give the following affirmative answer.

\begin{introthm}\label{introthm:maximal-converse}
Let $X$ be a discrete metric space of bounded geometry and $\H_0$ be 
a nonzero separable Hilbert space. Then the canonical homomorphism
\[
 \lambda_{X,\H_0}:C^*_{\max}(X;\H_0)
 \longrightarrow C^*(X;\H_0)
\]
is an isomorphism \emph{if and only if} $X$ has Property~A.
\end{introthm}

In the language of coarse groupoids, Theorem~\ref{introthm:maximal-converse} reads as follows.

\begin{introcor}
Let $X$ be a discrete metric space of bounded geometry. Then $X$ has Property~A \emph{if and only if} the coarse groupoid $G(X)$ has the weak containment property.
\end{introcor}

Theorem~\ref{introthm:maximal-converse} also answers the questions
concerning maximal--reduced equality in
\cite[Examples~1.15 and~1.16]{SW}. Both the box space of
$\mathrm{SL}(2,\mathbb Z)$ considered in \cite[Example~1.15]{SW}
and the box space of a free group constructed by Arzhantseva,
Guentner and \v{S}pakula \cite{AGS} fail Property~A. Consequently,
the canonical quotient is not injective in either the
uniform or the ordinary Roe setting. Since the latter box space coarsely embeds into Hilbert space, the quotient
maps induce isomorphisms in K-theory by \cite[Theorem~1.4]{SW}. Their kernels are therefore nonzero ideals
with trivial K-theory.

As an immediate consequence of Theorem~\ref{introthm:maximal-converse}, we obtain the following:

\begin{introcor}
Let $X$ be a discrete metric space of bounded geometry and $\H_0$ be a nonzero separable Hilbert space. Then the algebraic Roe algebra $\CC[X;\H_0]$ has a unique $C^*$-norm \emph{if and only if} $X$ has Property~A.
\end{introcor}

To see this, it suffices to show that the concrete operator norm on $\CC[X;\H_0]$ is the smallest $C^*$-norm. Indeed, for every finite-rank projection $P$ supported on a finite subset of $X$, the corner $P\CC[X;\H_0]P$ is a full matrix algebra, on which every $C^*$-norm agrees with the concrete operator norm. Since such projections converge strongly to the identity, it follows that the concrete operator norm is dominated by every $C^*$-norm on $\CC[X;\H_0]$. Note that this contrasts with the situation for group rings, where the equality of the full and reduced $C^*$-norms does not in general imply uniqueness of the $C^*$-norm on the algebraic group ring, even for amenable groups. Indeed, the algebraic group ring $\CC[\ZZ]$ does not even admit a smallest $C^*$-norm \cite[Example~6.6]{GMR18}.


For discrete groups, Theorem~\ref{introthm:maximal-converse}
resolves two questions in \cite[Section~7, Questions~(2) and (3)]{BEW}. Denote by 
$C^*_{\max}(G)$ the full group $C^*$-algebra of a discrete
group $G$, and by $\rtimes_{\mathrm{inj}}$ the maximal injective
crossed-product functor introduced in \cite{BEW}.

\begin{introcor}\label{introcor:LLP}
Let $G$ be a discrete group. Then:
\begin{enumerate}
\item If $C^*_{\max}(G)$ has the local lifting property (LLP),
then $G$ is exact.
\item The maximal injective crossed-product functor
$\rtimes_{\mathrm{inj}}$ agrees with the reduced
crossed-product functor $\rtimes_r$ \emph{if and only if} $G$ is exact.
\end{enumerate}
\end{introcor}

The proof is given in Section~\ref{sec:maximal}, together with a
restricted tensor-norm characterisation of exactness for countable
discrete groups.

\subsection{Outline of proofs}

The converses in Theorems~\ref{introthm:uniform-converse} and
\ref{introthm:maximal-converse} are proved by contraposition, using
a common set of geometric and operator-theoretic tools. In
Section~\ref{sec:preparations}, failure of Property~A is used to
obtain widely separated finite weak expander pieces with uniform
local spectral estimates, and we also reduce the two
converses from arbitrary nonzero coefficient Hilbert spaces to the scalar
case by means of rank-one corners.

Section~\ref{sec:spec balance} establishes a matrix discrepancy
estimate and a spectral balancing result, which are applied in two
different ways. For Theorem~\ref{introthm:uniform-converse},
spectral balancing is combined with the operators from weak expanders to
construct uniformly quasi-local row operators. Finite-dimensional
representation theory and an argument of Ozawa then produce a
quasi-local operator which stays at positive distance from every
finite-propagation operator. The construction can even be arranged to
produce a ghost operator, which is carried out in Section~\ref{sec:proof-quasi-local}.

For Theorem~\ref{introthm:maximal-converse}, the local spectral
estimates from weak expanders lead to finite graphs with uniformly bounded degree,
large mean degree and girth tending to infinity. The
matrix discrepancy estimate supplies edge signings for which the
associated adjacency operators have small concrete norm. By
contrast, signed distance kernels on large-girth graphs give a
lower bound for their maximal norm, using a local positivity
estimate and the sum-of-squares criterion from \cite{MOSS}. 
This proves a strict maximal-reduced norm gap, and hence
the failure of weak containment, when Property~A fails
(see Section~\ref{sec:maximal}).

Finally, Theorem~\ref{introthm:sq=Roe} is proved independently in
Section~\ref{sec:strong quasiloc}. The main step is a quantitative approximation theorem
based on the weak-$*$ Hahn-Banach separation theorem, combined with a Lipschitz transport
lemma. These provide commutator estimates for
strongly quasi-local operators and control their distance from
finite-propagation operators.

\subsection*{AI use statement}
The problems studied in this paper originate from the authors' long-standing research on Roe algebras, quasi-locality, and related questions in coarse geometry. The authors had investigated these problems extensively and developed substantial related results before using AI, providing the mathematical background and motivation for the present work.
Building on this research, the authors used OpenAI's ChatGPT (including GPT-6) to explore and develop several technical arguments. In particular, ChatGPT assisted with the spectral balancing and matrix discrepancy arguments for Theorem A, the Lipschitz transport lemma and Hahn-Banach duality argument for Theorem D, and parts of the random graph construction and signed distance kernel arguments for Theorem E. ChatGPT also assisted in refining the exposition. All mathematical arguments were carefully examined, corrected and refined by the authors, who take full responsibility for the results and the final manuscript. We also use Claude Opus 5.5 for proofreading and grammatical corrections.

\section{Preparations: weak expanders and reduction}\label{sec:preparations}

The converses in Theorems~\ref{introthm:uniform-converse} and
\ref{introthm:maximal-converse} share two preliminary ingredients:
weak expanders and a reduction to scalar coefficients. We record
them in this section.

\subsection{Weak expanders and local spectral gaps}

Recall that a metric space $(X,d)$ has \emph{Property~A} if for any $r,\varepsilon>0$ there are $s>0$ and $\{\xi_x\}_{x\in X}$ in $\ell^1(X)$ with $\|\xi_x\|_1=1$ such that $\supp(\xi_x) \subseteq B(x,s)$ for any $x\in X$, and $\|\xi_x-\xi_y\|_1<\varepsilon$ whenever $d(x,y) \leq r$.

Let us fix some notation. Given a metric space $(X,d)$ and $F\subseteq X$, denote by $G_r(F)$ the graph on $F$ where distinct
points are joined by an edge if their distance in $X$ is at most $r$.
For $A\subseteq F$, write $\partial_r^F A:=\{y\in F\setminus A:d(y,A)\le r\}$.
For $r \geq 0$, denote $N(r):=\sup_{x\in X}|B(x,r)|$. For $A\subseteq X$ and $r\geq0$, denote the $r$-neighbourhood by $\mathcal N_r(A):=\{x\in X:d(x,A)\leq r\}$.
Given a Hilbert space $\H_0$ and $T\in\B(\ell^2(X;\H_0))$, define its \emph{propagation} to be $\ppg(T):=\sup\{d(x,y):T_{x,y}\neq0\}$, with $\ppg(0):=0$. Set
\[
q_r(T):=
\sup\{\|\chi_A T\chi_B\|: A,B \subseteq X \text{ with } d(A,B) > r\}.
\]

Let us recall the notion of weak expanders from \cite[Definition~2.4]{Sako}. 
\begin{defn}\label{defn:weak expander}
A sequence $(Y_n)_{n\in\NN}$ of non-empty finite subspaces of $X$ is called a \emph{weak expander sequence} if there exist $r_0,h>0$ satisfying the following: for every $s>0$ there exists $N \in \NN$ such that for any $n\geq N$, we have
\begin{equation}\label{eq:weak-diameter}
 |\partial_{r_0}^{Y_n} A|\ge h|A|
 \quad \text{for any non-empty }  A\subseteq Y_n \text{ with } \diam(A) \le s.
\end{equation}
\end{defn}


The following result characterises the failure of Property~A by the existence of weak expanders, which follows from uniform local amenability (ULA). Indeed, Property~A implies ULA \cite[Lemma~2.7 and Proposition~3.2]{BNSWW}, and the converse follows from \cite{Elek,ET}. On the other hand, Sako characterised the failure of ULA by the existence of a weak expander sequence $\{Y_n\}_n$ in $X$ \cite[Proposition~6.2]{Sako} (see also \cite[Definition~2.6]{RW}). For convenience, we include a proof showing that the $Y_n$ can be chosen pairwise disjoint with $d_X(Y_n, Y_m) \to \infty$ as $n+m \to \infty$ for $n \neq m$.

\begin{prop}\label{prop:pieces}
Let $X$ be a discrete metric space of bounded geometry.
Then $X$ does not have Property~A \emph{if and only if} there exist constants $r_0,h>0$ and a sequence of mutually disjoint non-empty finite subsets $\{Y_n\}_{n\in \NN}$ of $X$ such that \eqref{eq:weak-diameter} holds for $Y_n$ with $s=n$ and 
\begin{equation}\label{EQ:sparse}
 d_X\Big(Y_n,\bigcup_{i<n}Y_i\Big)>n \quad \text{for} \quad n\ge 2.
\end{equation}
\end{prop}

\begin{proof}
If such a sequence $\{Y_n\}_n$ exists, \cite[Corollary~6.3]{Sako} shows that $X$ does not have Property~A. Hence, we focus on the converse.

Suppose that $X$ fails Property~A. Fix $o\in X$ and divide $X$ into the annuli
\[
 X_m:=\{x\in X:m^2\leq d_X(o,x)<(m+1)^2\} \quad \text{for} \quad m=0,1,2,\cdots.
\]
Since Property~A passes to finite unions, either $\bigsqcup_{m:\text{ even}} X_m$ or $\bigsqcup_{m: \text{ odd}} X_m$ fails Property~A, denoted by $Z$. Clearly, for each $r>0$ the components of the induced graph $G_r(Z)$ are finite.

Since $Z$ does not have Property~A, it follows from the definition that there exists $r_0>0$ such that the graph $G_{r_0}(Z)$ does not have Property~A in the sense of \cite{ET}. By \cite[Theorem~1]{ET} (see also \cite[Theorem~1]{Elek}), $G_{r_0}(Z)$ does not have ULA in the sense of \cite[Section~1.2(3)]{ET}. Equivalently (since $G_{r_0}(Z)$ has bounded valency), there exists $c>0$ such that for any $M\in\NN$ there is a non-empty finite set $F\subseteq Z$ satisfying: for any non-empty $A\subseteq F$ with $|A|\le M$, we have $|\partial^F_{r_0}A| \geq c|A|$. Using bounded geometry of $X$, we obtain a sequence $\{F_j\}_{j\in \NN}$ of non-empty finite subsets of $Z$ such that for any non-empty $A\subseteq F_j$ with $\diam(A) \le j$, we have $|\partial^{F_j}_{r_0}A| \geq c|A|$.

Fix a connected component $C$ of $G_{r_0}(Z)$, and let $j\ge\diam(C)$. If $A:=F_j\cap C$ were non-empty, then $\diam(A) \leq j$ and $\partial_{r_0}^{F_j} A=\varnothing$, since every point of $F_j$
within distance $r_0$ of $A$ belongs to $C$. This contradicts the
choice of $F_j$. Thus, we obtain $F_j \cap C = \emptyset$. 

Set $Y_1:=F_1$. Having chosen $Y_1,\cdots,Y_{n-1}$, take $j\ge n$ large enough so that $F_j$ has empty intersection with the closed $n$-neighbourhood in $Z$ of $\bigcup_{i=1}^{n-1} Y_i$. Setting $Y_n:=F_j$, we finish the construction. 
\end{proof}

The following spectral consequence of weak expanders is essentially contained in \cite[Proposition 3.1]{RW}. Here we sketch the proof for the reader’s convenience.

\begin{prop}\label{prop:weak expander spectral}
Let $X$ be a discrete metric space of bounded geometry, which does not have Property~A. Then there exists a sequence of subspaces $\{Y_n\}_n$ satisfying \eqref{EQ:sparse} and a sequence of positive contractions $\{S_n\}_n$ with $S_n \in \B(\ell^2(Y_n))$ such that the following hold:
\begin{enumerate}
 \item $S_ne_n=e_n$ for $e_n:=|Y_n|^{-1/2}\chi_{Y_n}$ and $n\in \NN$;
 \item $\sup_n q_s(S_n) \to 0$ as $s \to +\infty$;
 \item $\max_{x\in Y_n}\|S_n\delta_x\| \to 0$ as $n \to +\infty$.
\end{enumerate}
\end{prop}

\begin{proof}[Sketch of proof]
Let us fix $r_0,h>0$, a sequence $\{Y_n\}_{n\in \NN}$ as in Proposition~\ref{prop:pieces} and set $D:=N(r_0)$. For each $n\in \NN$, set $H_n:=G_{r_0}(Y_n)$, which has valency at most $D$. Set $\Delta_n$ to be the graph Laplacian, \emph{i.e.}, 
\[
(\Delta_n f)(x)=\sum_{y\sim_{H_n}x}\bigl(f(x)-f(y)\bigr) \quad \text{for} \quad f \in \ell^2(Y_n).
\]
Here $y\sim_{H_n}x$ means that $x,y$ are adjacent in $H_n$. Put $L_n:=\frac{\Delta_n}{2D}$ and $\mu:= \frac{h^2}{4D^2}$.  Clearly, $0\le L_n\le \Id_{\ell^2(Y_n)}$ since each vertex in $H_n$ has valency at most $D$. 

Next, we show that for every non-empty $B \subseteq Y_n$ with diameter at most $n$, the spectrum of the compression $\chi_B L_n \chi_B \in \B(\ell^2(B))$ lies in $[\mu,1]$. To see this, let $f\in\ell^2(Y_n)$ be a unit vector with $\diam(\supp f)\leq n$, and put $\varphi=|f|^2$. For $t\geq 0$, let $A_t:=\{x\in Y_n:\varphi(x)>t\}$. Each non-empty $A_t$ has diameter at most $n$, so Proposition~\ref{prop:pieces}
gives $|\partial_{r_0}^{Y_n}A_t|\geq h|A_t|$. Thus, 
\[
  \sum_{x\in Y_n}\sum_{y\sim_{H_n}x}
  |\chi_{A_t}(x)-\chi_{A_t}(y)|
  \geq 2|\partial_{r_0}^{Y_n}A_t|
  \geq 2h|A_t|.
\]
Integrating over $t$ gives
\[
  \sum_{x\in Y_n}\sum_{y\sim_{H_n}x}|\varphi(x)-\varphi(y)|
  \geq 2h\sum_{x\in Y_n}\varphi(x)=2h.
\]
This is inequality (3.1) in the proof of \cite[Lemma~3.2]{RW} with $c=2h$. Their estimate with $M=D$ gives
$\langle\Delta_n f,f\rangle\geq \frac{(2h)^2}{8D} = \frac{h^2}{2D}$ and hence, $\langle L_nf,f\rangle\ge\mu$.

Finally, choose a continuous function $\varphi:[0,1]\to[0,1]$ with $\varphi(0)=1$ and $\varphi|_{[\mu,1]}=0$. For each $n\in \NN$, set $S_n:=\varphi(L_n)$. Since $L_ne_n=0$ and $\varphi(0)=1$, we obtain condition (1). Let $Y=\bigcup_n Y_n$. The operator $L:=\mathrm{(SOT)}$-$\sum_{n\in \NN} L_n$ on $\ell^2(Y)$ has propagation at most $r_0$. Hence, $\varphi(L)=\mathrm{(SOT)}$-$\sum_{n\in \NN} S_n$ belongs to $C_u^*(Y)$ and is quasi-local. Since $q_s(S_n)\leq q_s(\varphi(L))$ for each $n$, we obtain condition (2). Condition (3) follows from the estimate in \cite[proof of Proposition~3.1]{RW}, with the help of the local spectral gap from the previous paragraph.
\end{proof}

\subsection{Reduction to the uniform case}\label{ssec:reductions}

We next reduce the general $\H_0$-case in both Theorems~\ref{introthm:uniform-converse} and \ref{introthm:maximal-converse} to the scalar case (\emph{i.e.}, $\H_0=\CC$) using the same rank-one corner map. Note that a related rank-one corner construction was used in \cite[Section~4.2]{SW} to establish the Morita equivalence between the maximal uniform Roe algebra and the maximal uniform algebra.

\begin{lem}\label{lem:rank-one-corner}
Let $\H_0$ be a nonzero separable Hilbert space with a unit vector
$e\in\H_0$, and $P\in\B(\H_0)$ be the orthogonal projection onto
$\CC e$. For $T\in\B(\ell^2(X))$, we define $j(T):=T\otimes P \in \B(\ell^2(X;\H_0))$. Then:
\begin{enumerate}
 \item For any $T\in\B(\ell^2(X))$, $j(T)$ is locally compact. Moreover, $T \in C^*_u(X)$ if and only if $j(T) \in C^*(X;\H_0)$, and $T\in C^*_{uq}(X)$ if and only if $j(T) \in C^*_q(X;\H_0)$;
 \item For any $T\in\CC_u[X]$, $\|j(T)\|=\|T\|$ and $\|j(T)\|_{\max}=\|T\|_{\max}$.
\end{enumerate}
\end{lem}

\begin{proof}
Let $J:\ell^2(X)\to\ell^2(X;\H_0)$ be the isometry given by $J\delta_x:=\delta_x\otimes e$ for each $x\in X$, where $\delta_x$ is the Dirac function at $x$. For each $T\in\B(\ell^2(X))$, we have $j(T)=JTJ^*$.

(1). For $T\in \B(\ell^2(X))$, clearly $T\otimes P \in \B(\ell^2(X;\H_0))$ is locally compact. Moreover, for $A,C\subseteq X$,
\[
 \big\|(\chi_A\otimes\Id_{\H_0})(T\otimes P)
       (\chi_C\otimes\Id_{\H_0})\big\|
 =\|\chi_AT\chi_C\|.
\]
Thus, $T\in C^*_{uq}(X)$ if and only if $j(T) \in C^*_q(X;\H_0)$. It is also clear that $T \in C^*_u(X)$ implies $j(T) \in C^*(X;\H_0)$. For the converse, choose locally compact finite-propagation operators $S_n$ on
$\ell^2(X)\otimes \H_0$ such that $\|S_n-j(T)\|\to0$. Since $J$ has propagation zero, each $J^*S_nJ$ has finite propagation. Moreover, $J^* j(T) J=T$ implies that
\[
 \|J^*S_nJ-T\|=\|J^*(S_n-j(T))J\|
 \leq\|S_n-j(T)\|\to 0 \quad \text{as} \quad n \to +\infty.
\]
Hence, we obtain that $T\in C_u^*(X)$.

(2). Let $T\in\CC_u[X]$. Clearly $\|j(T)\|=\|T\|$, so it remains to consider the maximal norm. The map
$j:\CC_u[X]\to \CC[X;\H_0]$ is a $*$-homomorphism, so composing representations with $j$ gives $\|j(T)\|_{\max} \leq \|T\|_{\max}$. 
For the converse, denote $Q:=JJ^*=\Id_{\ell^2(X)}\otimes P \in \CC[X;\H_0]$. Define the map
\[
 E:\CC[X;\H_0]\longrightarrow\CC_u[X]
 \quad \text{by} \quad  E(B):=J^*BJ.
\]

Claim: for each $B \in \CC[X;\H_0]$, $E(B^*B)$ is a finite sum of squares in $\CC_u[X]$.
Fix $B\in \CC[X;\H_0]$ with $r:=\ppg(B)$ and put $N:=N(r)$.
For each $z\in X$, let
\[
 W_z:=\operatorname{span}\{B_{z,x}e:x\in B(z,r)\}\subseteq\H_0
 \quad \text{and} \quad k_z:=\dim W_z\le N.
\]
Choose an orthonormal basis $u_{z,1},\ldots,u_{z,k_z}$ of $W_z$,
and set $u_{z,i}=0$ for $k_z<i\le N$. For $1\le i\le N$, define $c_i$ by 
\[
 (c_i)_{z,x}:=\langle B_{z,x}e,u_{z,i}\rangle
 \quad \text{for} \quad z,x\in X.
\]
Each $c_i$ has propagation at most $r$ and bounded entries, hence
$c_i\in\CC_u[X]$. Moreover, for $x,y\in X$ we have
\[
 E(B^*B)_{x,y}
 =\langle BJ\delta_y,BJ\delta_x\rangle
   =\sum_{z\in X}\langle B_{z,y}e,B_{z,x}e\rangle
 =\sum_{i=1}^{N}\sum_{z\in X}
     \overline{(c_i)_{z,x}}(c_i)_{z,y}
   =\sum_{i=1}^{N}(c_i^*c_i)_{x,y}.
\]
Thus
\begin{equation}\label{eq:corner-squares}
 E(B^*B)=\sum_{i=1}^{N}c_i^*c_i.
\end{equation}

Let $\omega$ be a state on $C^*_{u,\max}(X)$, and define
$F:\CC[X;\H_0]\to\CC$ by $F(B):=\omega(E(B))$.
By~\eqref{eq:corner-squares}, $F(B^*B)\ge0$ for every
$B\in\CC[X;\H_0]$. We verify that left multiplication in the GNS
construction for $F$ is bounded. Fix $S\in\CC[X;\H_0]$ of propagation
at most $r'$, and put $N':=N(r')$. By \cite[Lemma~2]{Oza}, we can write $S=D_1+\cdots+D_{2N'}$, where each $D_i\in\CC[X;\H_0]$ has at most one nonzero matrix entry
in each row and column. In particular, $D_i^*D_i$ is diagonal and
$\|D_i\|\le\|S\|$. Fix $c>\|S\|$ and set
\[
 R_i:=(c^2\Id_{\ell^2(X;\H_0)}-D_i^*D_i)^{1/2}
 \quad \text{for} \quad 1\le i\le 2N'.
\]
Since $D_i^*D_i$ is diagonal and has compact diagonal entries, 
$R_i-c\Id_{\ell^2(X;\H_0)}\in\CC[X;\H_0]$.
Thus $R_iB\in\CC[X;\H_0]$ for every $B\in\CC[X;\H_0]$, and the
positivity of $F$ gives
\[
0\le F((R_iB)^*(R_iB))=c^2F(B^*B)-F((D_iB)^*(D_iB)).
\]
Writing $[B]$ for the class of $B$ in the GNS pre-Hilbert space,
we obtain
\[
 \|[D_iB]\|^2=F((D_iB)^*(D_iB))\le c^2F(B^*B)
 =c^2\|[B]\|^2.
\]
Hence, 
\[
 \|[SB]\|\le\sum_{i=1}^{2N'}\|[D_iB]\|
 \le 2N'c\|[B]\|\quad \text{for each} \quad B\in\CC[X;\H_0].
\]
Therefore, the GNS construction yields a $*$-representation
$\pi_F:\CC[X;\H_0]\to\B(\H_F)$ on the corresponding Hilbert space
$\H_F$.

Since $E(Q)=\Id_{\ell^2(X)}$ and $Q^*Q=Q$, the vector
$\xi_F:=[Q]\in\H_F$ satisfies
\[
 \|\xi_F\|^2=F(Q)=\omega(\Id_{\ell^2(X)})=1.
\]
For $T\in\CC_u[X]$, we have $j(T)Q=j(T)$ and
$E(j(T)^*j(T))=T^*T$. Consequently,
\[
 \omega(T^*T)
 =F(j(T)^*j(T))
   =\|\pi_F(j(T))\xi_F\|^2
   \le\|\pi_F(j(T))\|^2
   \le\|j(T)\|_{\max}^2.
\]
Taking the supremum over all states $\omega$ gives
$\|T\|_{\max}\le\|j(T)\|_{\max}$.
\end{proof}

Consequently, we obtain the following:

\begin{cor}\label{cor:scalar-reductions}
Let $\H_0$ be a nonzero separable Hilbert space. Then:
\begin{enumerate}
\item If $C^*(X;\H_0)=C_q^*(X;\H_0)$, then
      $C_u^*(X)=C_{uq}^*(X)$.
\item If $\lambda_{X,\H_0}:C_{\max}^*(X;\H_0)\to C^*(X;\H_0)$
      is injective, so is 
      $\lambda_X:C^*_{u,\max}(X)\to C_u^*(X)$.
\end{enumerate}
\end{cor}

Thus, for the converses of Theorems~\ref{introthm:uniform-converse}
and~\ref{introthm:maximal-converse}, it suffices to work in the scalar case $\H_0=\CC$.

%

\section{Spectral balancing and matrix discrepancy}\label{sec:spec balance}

We establish a spectral balancing result using matrix discrepancy.
The partition estimate \eqref{eq:range-cells} will be used in Section~\ref{sec:proof-quasi-local},
while the underlying signing lemma (Lemma~\ref{lem:range-signing}) will also be used to control signed
adjacency matrices in Section~\ref{sec:maximal}. 

We use the following notation. For a vector \(u\) in a complex Hilbert space \(\H\), denote by \(\theta_u\in\B(\H)\) the positive rank-one operator defined by
\begin{equation}\label{EQ:theta operator}
\theta_u(w):=\langle w,u\rangle u \quad \text{for} \quad w\in \H.
\end{equation}
For a trace-class operator $A$ on $\H$, denote by $\|A\|_1:=\Tr(|A|)$ its trace norm.


The following proposition provides a partition that satisfies the same estimate for every self-adjoint test operator.

\begin{prop}\label{prop:range-split}
Given $k\in\NN$ with $k\geq 1$ and $\tau,\gamma>0$, there exists
$\ell>0$ satisfying the following property. Let \(Y\) be a finite set and \(\{u_x\}_{x\in Y}\) a family of vectors in a finite-dimensional Hilbert space \(\H\) such that
\[
 \sum_{x\in Y}\theta_{u_x}\leq\Id_{\H},
 \quad\text{and}\quad
 \|u_x\|^2\leq\ell\quad\text{for any} \quad x\in Y.
\]
Then there exist pairwise disjoint subsets $C_\sigma\subseteq Y$,
indexed by $\sigma\in\{-1,1\}^k$, such that
$B:=Y\setminus\bigcup_\sigma C_\sigma$ satisfies $|B|\leq\tau|Y|$.
Moreover, for each $\sigma$ and each self-adjoint operator $A$ on $\H$,
we have
\begin{equation}\label{eq:range-cells}
 \sum_{x\in C_\sigma}|\langle Au_x,u_x\rangle|
 \leq 2^{-k}\sum_{x\in Y}|\langle Au_x,u_x\rangle|+\gamma\|A\|_1.
\end{equation}
\end{prop}

To prove Proposition~\ref{prop:range-split}, we introduce a quantitative matrix discrepancy estimate. The lower bound on the number of sign vectors will be used to obtain a partition with a controlled discarded set.

\begin{lem}\label{lem:range-signing}
Let $E$ be a non-empty finite set, $d\geq1$ an integer and $\ell>0$.
Suppose that $\{u_x\}_{x\in E}$ is a family of vectors in a finite-dimensional
Hilbert space $\H$ satisfying
$\sum_{x\in E}\theta_{u_x}\leq\Id_{\H}$ and
$\|u_x\|^2\leq\ell$ for each $x\in E$.
Put $m:=|E|$ and $p:=\lceil m/d\rceil$. Then there are at least $2^{m-p}$
elements $s=(s_x)_{x\in E}\in\{-1,1\}^E$ (called \emph{sign vectors}) satisfying
\begin{equation}\label{EQ:sign vector}
 \Big\|\sum_{x\in E}s_x\theta_{u_x}\Big\|
 \leq8d^2\sqrt\ell.
\end{equation}
\end{lem}

\begin{proof}
Partition $E=P_1\sqcup\cdots\sqcup P_p$, where
$P_j=\{x_{j,1},\ldots,x_{j,m_j}\}$ and $1\leq m_j\leq d$.
Fix a sign vector $b\in\{-1,1\}^E$ and let
$\widetilde\H:=\bigoplus_{i=1}^d\bigoplus_{\varepsilon\in\{-1,1\}}\H$ be the $(2d)$-copy of $\H$.
For $1\leq j\leq p$, define $v_j\in\widetilde\H$ by
\[
 (v_j)_{i,\varepsilon}:=
 \begin{cases}
 d^{-1/2}u_{x_{j,i}},& i\leq m_j\ \text{and}\
 \varepsilon=b_{x_{j,i}},\\
 0,&\text{otherwise},
 \end{cases}
 \quad \text{for} \quad  1\leq i\leq d \quad \text{and} \quad \varepsilon\in\{-1,1\}.
\]
Then $\|v_j\|^2=d^{-1}\sum_{i=1}^{m_j}\|u_{x_{j,i}}\|^2\leq\ell$.
For $z=(z_{i,\varepsilon})\in\widetilde\H$, we have
\begin{align*}
\sum_{j=1}^p|\langle z,v_j\rangle|^2
&=\frac1d\sum_{j=1}^p
\left|\sum_{i=1}^{m_j}
\langle z_{i,b_{x_{j,i}}},u_{x_{j,i}}\rangle
\right|^2\leq\sum_{j=1}^p\sum_{i=1}^{m_j}
|\langle z_{i,b_{x_{j,i}}},u_{x_{j,i}}\rangle|^2\\
&\leq\sum_{i=1}^d\sum_{\varepsilon\in\{-1,1\}}
\sum_{x\in E}|\langle z_{i,\varepsilon},u_x\rangle|^2\leq\sum_{i=1}^d\sum_{\varepsilon\in\{-1,1\}}
\|z_{i,\varepsilon}\|^2
=\|z\|^2,
\end{align*}
where we use the Cauchy--Schwarz inequality in the first inequality and the assumption $\sum_{x\in E}\theta_{u_x}\leq\Id_{\H}$ in the last inequality. 
Consequently, for $A_j:=\theta_{v_j}$ we have
$\sum_jA_j\leq\Id_{\widetilde\H}$ and
$\sum_jA_j^2\leq\ell\Id_{\widetilde\H}$.

By \cite[Corollary~1.6]{KLS}, there exist
$t_1,\ldots,t_p\in\{-1,1\}$ such that
\[
\|T\|\leq4\sqrt\ell \quad \text{for} \quad T:=\sum_{j=1}^pt_jA_j.
\]
For $1\leq i\leq d$ and $\varepsilon\in\{-1,1\}$,
let $J_{i,\varepsilon}:\H\to\widetilde\H$ be the canonical
isometric inclusion into the $(i,\varepsilon)$-summand. Then for \(1\le i\le m_j\), direct calculations show that
\[
d\left(
J_{i,1}^*A_jJ_{i,1}
-
J_{i,-1}^*A_jJ_{i,-1}
\right)
=b_{x_{j,i}}\theta_{u_{x_{j,i}}}.
\]
Then we have
\[
 \sum_{j=1}^p\sum_{i=1}^{m_j}
 t_jb_{x_{j,i}}\theta_{u_{x_{j,i}}}
 =d\sum_{i=1}^d
 \bigl(J_{i,1}^*TJ_{i,1}-J_{i,-1}^*TJ_{i,-1}\bigr),
\]
whose norm is at most $2d^2\|T\|\leq8d^2\sqrt\ell$.

In conclusion, we showed that for each $b\in\{-1,1\}^E$, there exist signs $t_j$ for $j=1,\cdots, p$ such that
$s_{x_{j,i}}:=t_jb_{x_{j,i}}$ satisfies \eqref{EQ:sign vector}. For $b,c\in\{-1,1\}^E$, write $b\sim c$ if
there exist $t_1,\cdots,t_p\in\{-1,1\}$ such that
$c_x=t_jb_x$ for all $x\in P_j$ and $1\leq j\leq p$.
Clearly, each equivalence class has exactly
$2^p$ elements. Since there are
$2^{m-p}$ equivalence classes and 
each equivalence class contains a sign vector
satisfying \eqref{EQ:sign vector}, at least
$2^{m-p}$ sign vectors satisfy \eqref{EQ:sign vector}.
\end{proof}

Now we use Lemma~\ref{lem:range-signing} to obtain a two-set decomposition and then iterate it to prove Proposition~\ref{prop:range-split}.

\begin{proof}[Proof of Proposition~\ref{prop:range-split}]
Choose an integer $d\geq1$ and $0<\ell\leq1$ such that
\[
 k\alpha\leq\tau \quad \text{for} \quad 
 \alpha:=2\sqrt{\frac{\log2}{d}},\quad \text{and} \quad 
 \eta:=8d^2\sqrt\ell\leq\gamma.
\]
For a self-adjoint $A\in\B(\H)$ and $F\subseteq Y$, denote 
$w_A(F):=\sum_{x\in F}|\langle Au_x,u_x\rangle|$. 

Firstly, we prove that each non-empty $E\subseteq Y$ admits
disjoint subsets $E_+,E_-\subseteq E$ satisfying $|E\setminus(E_+\cup E_-)|\leq\alpha|E|$ and the following:
\begin{equation}\label{eq:one-splitting}
 w_A(E_\pm)\leq\tfrac12w_A(E)+\tfrac\eta2\|A\|_1
 \quad \text{for any} \quad  A=A^*\in\B(\H).
\end{equation}
To see this, let us fix a non-empty $E\subseteq Y$ and put $m:=|E|$.
If $m<d$, choose any partition $E=E_+\sqcup E_-$. Then we have
\begin{align*}
w_A(E_\pm) &\le w_A(E)=\frac12w_A(E)+\frac12w_A(E)\le\frac12w_A(E)+\frac12\sum_{x\in E}\|A\| \cdot \|u_x\|^2\\
& \leq \frac12w_A(E) + \frac{m\ell}{2}\|A\|_1 \leq \frac12w_A(E) + \frac{\eta}{2}\|A\|_1.
\end{align*}
Now assume that $m\geq d$, and set $p:=\lceil m/d\rceil$. By Lemma~\ref{lem:range-signing}, there are at least $2^{m-p}$ sign vectors $s = (s_x)_{x\in E}\in\{-1,1\}^E$ such that for any $A=A^*$ in $\B(\H)$, 
\[
\left|\sum_{x\in E}s_x\langle Au_x,u_x\rangle\right| =\left|\Tr\left(
A\sum_{x\in E}s_x\theta_{u_x}
\right)\right|\le\left\|\sum_{x\in E}s_x\theta_{u_x}\right\|\cdot\|A\|_1\le\eta\|A\|_1,
\]
where the first equality follows from $\Tr(A\theta_{u_x})=\langle Au_x,u_x\rangle$ by a direct calculation. 

Claim: there are at least $2^{m-p}$ sign vectors $s= (s_x)_{x\in E}\in\{-1,1\}^E$ such that
\begin{equation}\label{EQ:Claim admissible signs}
\sum_{x\in E}s_x|\langle Au_x,u_x\rangle|\le\eta\|A\|_1 \quad \text{for any} \quad A=A^* \quad \text{in} \quad \B(\H).
\end{equation}
To see this, write $E=\{x_1,\ldots,x_m\}$. We replace
$\langle Au_{x_j},u_{x_j}\rangle$ by $|\langle Au_{x_j},u_{x_j}\rangle|$, one index at a time. Fix $1\leq j\leq m$ and the signs
$s_{x_i}$ for $i\neq j$. For $A=A^*$, put
\[
c_A:=
\sum_{i=1}^{j-1}s_{x_i}|\langle Au_{x_i},u_{x_i}\rangle|
+\sum_{i=j+1}^m s_{x_i}\langle Au_{x_i},u_{x_i}\rangle.
\]
If exactly one choice $s_{x_j}=\varepsilon\in\{-1,1\}$ satisfies 
\begin{equation}\label{EQ:admissible}
c_A+\varepsilon \langle Au_{x_j},u_{x_j}\rangle\leq\eta\|A\|_1 \quad \text{for any} \quad A=A^*\in \B(\H),
\end{equation}
then $c_A-|\langle Au_{x_j},u_{x_j}\rangle|\leq c_A+\varepsilon \langle Au_{x_j},u_{x_j}\rangle \leq\eta\|A\|_1$ and we replace $s_{x_j}$ by $-1$. If both choices $s_{x_j}=\pm1$ satisfy \eqref{EQ:admissible}, then for any self-adjoint $A \in \B(\H)$ we have
\[
c_A\pm|\langle Au_{x_j},u_{x_j}\rangle|\leq\max\{c_A+\langle Au_{x_j},u_{x_j}\rangle,c_A-\langle Au_{x_j},u_{x_j}\rangle\}\leq\eta\|A\|_1,
\]
so we keep both $s_{x_j}=\pm1$. Thus, for each fixed choice of the other signs, the number of chosen $s_{x_j}$ cannot decrease, which proves the claim by induction on $j$.

Denote the set $\G:=\{s\in\{-1,1\}^E: \eqref{EQ:Claim admissible signs} \text{ holds for }s\}$, and equip $\{-1,1\}^E$ with the Hamming metric $d_H(s,t):=|\{x\in E:s_x\ne t_x\}|$. Now we use concentration of measure on the
discrete cube to show that $\diam_{d_H}(\G) \geq (1-\alpha)m$. 
Define
\[
\Phi(r):=d_H(r,\mathcal G) \quad \text{for} \quad r\in\{-1,1\}^E, \quad \text{and} \quad \overline{\Phi}:=\frac{1}{2^m}\sum_{r\in\{-1,1\}^E}\Phi(r).
\]
Since $\Phi$ is $1$-Lipschitz, McDiarmid's concentration inequality on the discrete cube (\cite[Theorem~6.7]{McD}) shows that for any $t \geq 0$, we have
\[
\frac{1}{2^m}
\left|\left\{r\in\{-1,1\}^E:
\Phi(r)\leq\overline{\Phi}-t\right\}\right|
\leq \exp\left(-\frac{2t^2}{m}\right).
\]
Taking $t=\overline{\Phi}$, we obtain
\[
2^{-p}
\leq\frac{|\mathcal G|}{2^m}
=\frac{|\{r\in\{-1,1\}^E:\Phi(r)=0\}|}{2^m}
\leq\exp\left(-\frac{2\overline{\Phi}^2}{m}\right),
\]
which implies that $\overline{\Phi}\leq\sqrt{\frac{mp\log 2}{2}}$.
Since $r\mapsto-r$ is a bijection on $\{-1,1\}^E$, we have
\[
\frac1{2^m}\sum_{r\in\{-1,1\}^E}
\bigl(\Phi(r)+\Phi(-r)\bigr)
=2\overline{\Phi}.
\]
Hence, there exists $r\in\{-1,1\}^E$ such that $\Phi(r)+\Phi(-r)\leq2\overline{\Phi}$. Choose $\hat{s},\hat{t}\in\mathcal G$ satisfying $d_H(r,\hat{s})=\Phi(r)$ and $d_H(-r,\hat{t})=\Phi(-r)$. Then
\begin{align*}
m-d_H(\hat{s},\hat{t}) &= d_H(\hat{s},-\hat{t}) \leq d_H(\hat{s},r)+d_H(r,-\hat{t}) = \Phi(r)+\Phi(-r) \leq 2\overline\Phi \\
&\leq \sqrt{2mp\log2} \leq 2m\sqrt{\frac{\log2}{d}} =\alpha m,
\end{align*}
where the last inequality follows from $p=\lceil m/d\rceil\leq2m/d$ since $m\geq d$.

For any $s\in\G$, denote $E_s:=\{x\in E: s_x=1\}$.
Set $E_+:=E_{\hat{s}}\setminus E_{\hat{t}}$ and $E_-:=E_{\hat{t}}\setminus E_{\hat{s}}$. Then $E_+$ and $E_-$ are disjoint, and we have
\[
|E\setminus(E_+\cup E_-)|=m-d_H(\hat{s},\hat{t})\leq\alpha m.
\]
Moreover, for any $s\in \G$ we have 
\[
\sum_{x\in E}s_x|\langle Au_x,u_x\rangle| =w_A(E_s)-w_A(E\setminus E_s)=2w_A(E_s)-w_A(E).
\]
Hence, $w_A(E_s)\le\frac12w_A(E)+\frac{\eta}{2}\|A\|_1$. The same estimate holds for every subset of $E_s$ and in particular for $E_+$ and $E_-$. This proves \eqref{eq:one-splitting}.

Finally, we apply the decomposition above recursively $k$ times, starting with $E=Y$. At the first step, we obtain
two disjoint subsets $E_+,E_-\subseteq Y$.
At each subsequent step, we apply the decomposition above to
every non-empty set obtained at the preceding step,
and replace an empty set by two empty sets. 
After $k$ steps, we obtain pairwise disjoint
sets $C_\sigma$, indexed by $\sigma\in\{-1,1\}^k$,
such that
\[
\left|Y\setminus
\bigcup_{\sigma\in\{-1,1\}^k}C_\sigma\right|
\leq k\alpha|Y|\leq\tau|Y|.
\]
Moreover, we have
\[
w_A(C_\sigma) \leq 2^{-k}w_A(Y)
+\frac{\eta}{2}\sum_{j=0}^{k-1}2^{-j}\|A\|_1=2^{-k}w_A(Y)
+\eta(1-2^{-k})\|A\|_1\leq 2^{-k}w_A(Y)+\gamma\|A\|_1
\]
for every $\sigma\in\{-1,1\}^k$ and 
self-adjoint $A\in\B(\H)$.
This completes the proof.
\end{proof}

\section{Quasi-local characterisation}
\label{sec:proof-quasi-local}

In this section, we prove Theorem~\ref{introthm:uniform-converse}. Our main ingredients are the operators from weak expanders (Proposition~\ref{prop:weak expander spectral}) and the spectral balancing result (Proposition~\ref{prop:range-split}). Let us divide the proof into several parts.

\subsection{Quasi-local estimates}

We first derive a quasi-locality estimate for row operators built
from the partitions in Section~\ref{sec:spec balance}, which will be applied to the operators from weak expanders in the next subsection. 

Let $Y$ be a non-empty finite set, $k\geq1$ be an integer and $d:=2^k$. Fix pairwise disjoint subsets $C_\sigma\subseteq Y$, indexed by $\sigma\in\{-1,1\}^k$, and set
\[
G:=\bigcup_{\sigma\in\{-1,1\}^k}C_\sigma \quad \text{and} \quad B:=Y\setminus G.
\]
Also fix a contraction $H\in\B(\ell^2(Y))$. We have the following:

\begin{prop}\label{prop:constructing quasi-local oeprators}
Suppose that there is $\gamma\geq0$ such that for any $\sigma\in\{-1,1\}^k$ and any self-adjoint $A\in\B(\ell^2(Y))$, we have
\begin{equation}\label{eq:abstract-balancing}
 \sum_{x\in C_\sigma}|\langle AH^*\delta_x,H^*\delta_x\rangle|
 \leq d^{-1}\sum_{x\in Y}|\langle AH^*\delta_x,H^*\delta_x\rangle|
       +\gamma\|A\|_1.
\end{equation}
For $j=1,\cdots,k$, define 
\[
V_j:=\sum_{\sigma\in\{-1,1\}^k}
\sigma_j\chi_{C_\sigma}H \in \B(\ell^2(Y)),
\]
and the row operator
\[
 V:\CC^k\otimes\ell^2(Y)\longrightarrow\ell^2(Y)\quad \text{by} \quad 
 V(\xi_1,\ldots,\xi_k):=\sum_{j=1}^kV_j\xi_j.
\]
Then for $E,D\subseteq Y$ with $h:=\|\chi_EH\chi_D\|$, we have
\begin{equation}\label{eq:abstract-row-bound}
 \|\chi_EV(\Id_{\CC^k}\otimes\chi_D)\|^2
 \leq h^2+\sqrt2\,h+2\varepsilon+\sqrt{2\varepsilon} \quad \text{for} \quad \varepsilon:=d\|\chi_BH\|^2+d^2\gamma.
\end{equation}
\end{prop}


For the proof, we first extend the balancing estimate from each $C_\sigma$ to arbitrary subsets $E_\sigma\subseteq C_\sigma$, using the Hahn-Banach theorem.

\begin{lem}\label{lem:range-compressions}
For every $\sigma\in\{-1,1\}^k$, we have $\|dH^*\chi_{C_\sigma}H-H^*H\|\leq\varepsilon$.
Moreover, for every $E_\sigma\subseteq C_\sigma$ there exists $f_\sigma:Y\to[0,1]$ such that
\begin{equation}\label{eq:range-subcell}
 \|dH^*\chi_{E_\sigma}H-H^*f_\sigma H\|\leq\varepsilon.
\end{equation}
\end{lem}

\begin{proof}
For $\xi\in\ell^2(Y)$, applying \eqref{eq:abstract-balancing} to $A=\theta_{\xi}$, we obtain
\[
 \|\chi_{C_\sigma}H\xi\|^2
 \leq d^{-1}\|H\xi\|^2+\gamma\|\xi\|^2
 \quad \text{for each} \quad \sigma\in\{-1,1\}^k.
\]
This also implies the following:
\begin{align*}
 \|\chi_{C_\sigma}H\xi\|^2
 &=\|H\xi\|^2-\|\chi_BH\xi\|^2
   -\sum_{\rho\ne\sigma}\|\chi_{C_\rho}H\xi\|^2\\
 &\geq \|H\xi\|^2-\|\chi_BH\xi\|^2 -(d-1)\left(d^{-1}\|H\xi\|^2+\gamma\|\xi\|^2\right)\\
 & = d^{-1}\|H\xi\|^2-\|\chi_BH\xi\|^2-(d-1)\gamma\|\xi\|^2.
\end{align*}
Combining them together, we obtain $\left|d\|\chi_{C_\sigma}H\xi\|^2-\|H\xi\|^2\right|\le\varepsilon\|\xi\|^2$. Since the operator \(dH^*\chi_{C_\sigma}H-H^*H\) is self-adjoint, the first statement holds.

Fix $E_\sigma\subseteq C_\sigma$ and suppose that
the second conclusion fails for $E_\sigma$. Note that the set $\{H^*fH ~|~ f:Y\to[0,1]\}$ is norm-compact and convex. Hence by the Hahn-Banach theorem, there is a self-adjoint $A\in\B(\ell^2(Y))$ with $\|A\|_1=1$
such that
\begin{equation}\label{EQ:separation}
\Tr(A \cdot dH^*\chi_{E_\sigma}H)
>\max_{f:Y \to [0,1]}\Tr(A \cdot H^*fH)
+\varepsilon.
\end{equation}
For $x\in Y$, set $a_x:=\langle AH^*\delta_x,H^*\delta_x\rangle$ and $(a_x)_+:=\max\{a_x,0\}$. Since
\[
\max_{f:Y \to [0,1]}\Tr(A \cdot H^*fH)
=\max_{f:Y \to [0,1]}\sum_{x\in Y}f(x)a_x
=\sum_{x\in Y}(a_x)_+,
\]
inequality \eqref{EQ:separation} becomes
\begin{equation}\label{eq:range-cell-separation}
d\sum_{x\in E_\sigma}a_x
>\sum_{x\in Y}(a_x)_++\varepsilon.
\end{equation}

On the other hand, \eqref{eq:abstract-balancing} implies $d\sum_{x\in C_\sigma}|a_x|\leq\sum_{x\in Y}|a_x|+d\gamma$. Also the first statement shows that
\[
d\sum_{x\in C_\sigma}a_x=\Tr(A \cdot dH^*\chi_{C_\sigma}H)\leq\operatorname{Tr}(AH^*H)+\varepsilon=\sum_{x\in Y}a_x+\varepsilon.
\]
Adding them together and using $2(a_x)_+=|a_x|+a_x$, we obtain
\[
d\sum_{x\in C_\sigma}(a_x)_+\leq\sum_{x\in Y}(a_x)_+
+\frac{d\gamma+\varepsilon}{2}\leq\sum_{x\in Y}(a_x)_++\varepsilon.
\]
Hence,
\[
d\sum_{x\in E_\sigma}a_x \leq d\sum_{x\in E_\sigma}(a_x)_+\leq d\sum_{x\in C_\sigma}(a_x)_+\leq\sum_{x\in Y}(a_x)_++\varepsilon,
\]
contradicting \eqref{eq:range-cell-separation}. This completes the proof.
\end{proof}

%


Now we prove Proposition~\ref{prop:constructing quasi-local oeprators}.

\begin{proof}[Proof of Proposition~\ref{prop:constructing quasi-local oeprators}]
Fix $\xi=(\xi_1,\ldots,\xi_k)$ with $\xi_j\in\ell^2(D)$ and $\sum_{j=1}^k\|\xi_j\|^2=1$.
For $\sigma\in\{-1,1\}^k$, set
\[
 \xi_\sigma:=\sum_{j=1}^k\sigma_j\xi_j
 \quad \text{and} \quad E_\sigma:=E\cap C_\sigma.
\]
Choose $f_\sigma:Y\to[0,1]$ satisfying \eqref{eq:range-subcell} for $E_\sigma$ by Lemma~\ref{lem:range-compressions}. Then we have
\begin{align*}
\|\chi_E V\xi\|^2
&=\sum_{\sigma\in\{-1,1\}^k}
  \|\chi_{E_\sigma}H\xi_\sigma\|^2=\frac1d\sum_\sigma
  \langle dH^*\chi_{E_\sigma}H\xi_\sigma,
  \xi_\sigma\rangle\\
  &\leq\frac1d\sum_\sigma
  \langle H^*f_\sigma H\xi_\sigma,
  \xi_\sigma\rangle
  +\frac{\varepsilon}{d}
  \sum_\sigma\|\xi_\sigma\|^2\\
&=\frac1d\sum_\sigma\sum_{x\in Y}
  f_\sigma(x)
  \left|\sum_{j=1}^k
  \sigma_j \cdot (H\xi_j)(x)\right|^2
  +\varepsilon,
 \end{align*}
where the last equality follows from $d^{-1}\sum_\sigma\|\xi_\sigma\|^2 =\sum_{j=1}^k\|\xi_j\|^2=1$.

Write $a(x):=\sum_{j=1}^k|(H\xi_j)(x)|^2$ and
$p(x):=d^{-1}\sum_\sigma f_\sigma(x)$.
Put $m:=\sum_{x\in Y}p(x)a(x)$. Applying
\eqref{eq:range-subcell} to each $\xi_j$ and summing over $\sigma$, we obtain
\[
 m=\frac1d\sum_\sigma\sum_{j=1}^k\sum_{x\in Y}
      f_\sigma(x)|(H\xi_j)(x)|^2\leq\sum_\sigma\sum_{j=1}^k\|\chi_{E_\sigma}H\xi_j\|^2
      +\varepsilon
   \leq h^2+\varepsilon.
\]
The last inequality follows since the $E_\sigma$ are disjoint
subsets of $E$ and $\xi_j\in\ell^2(D)$.
For $1\leq i<j\leq k$, put
$c_{ij}(x):=d^{-1}\sum_\sigma f_\sigma(x)\sigma_i\sigma_j$.
Then for $x\in Y$, we have
\[
 \frac1d\sum_\sigma f_\sigma(x)
 \left|\sum_{j=1}^k\sigma_j \cdot (H\xi_j)(x)\right|^2
 =p(x)a(x)+2\operatorname{Re}\sum_{i<j}c_{ij}(x)  \cdot
      (H\xi_i)(x)\overline{(H\xi_j)(x)}.
\]
The functions $\sigma\mapsto\sigma_i\sigma_j$ for $i<j$ are
orthonormal in $\ell^2(\{-1,1\}^k)$ with respect to the normalised counting
measure. By Bessel's inequality,
\[
 \sum_{i<j}|c_{ij}(x)|^2
 \leq\frac1d\sum_\sigma f_\sigma(x)^2\leq p(x).
\]
Consequently, the Cauchy--Schwarz inequality implies
\begin{align*}
2\left|\sum_{i<j}c_{ij}(x) \cdot
(H\xi_i)(x)\overline{(H\xi_j)(x)}\right| & \leq 2\sqrt{p(x)}
\left(\sum_{i<j}|(H\xi_i)(x)|^2 \cdot |(H\xi_j)(x)|^2\right)^{1/2}\\
&\leq
2\sqrt{p(x)} \cdot
\frac{1}{\sqrt 2}
\left(\sum_{i=1}^k|(H\xi_i)(x)|^2\right)=\sqrt{2p(x)}\,a(x).
\end{align*}
Since $\sum_{x\in Y}a(x)=\sum_{j=1}^k\|H\xi_j\|^2\leq1$,
we have
\[
 \sum_{x\in Y}\sqrt{p(x)}\,a(x)
 \leq\left(\sum_{x\in Y}p(x)a(x)\right)^{1/2}
      \left(\sum_{x\in Y}a(x)\right)^{1/2}\leq\sqrt m.
\]
Combining the above together, finally we obtain
\[
 \|\chi_EV\xi\|^2
 \leq m+\sqrt{2m}+\varepsilon\leq h^2+\sqrt2\,h+2\varepsilon+\sqrt{2\varepsilon}.
\]
Taking the supremum over $\xi$, we finish the proof.
\end{proof}

\subsection{Proof of Theorem~\ref{introthm:uniform-converse}}
This whole subsection is devoted to the proof of Theorem~\ref{introthm:uniform-converse}. 
Let us fix a nonzero separable Hilbert space $\H_0$. If $X$ has Property~A, \cite[Theorem~3.3]{SZ20} shows that $C^*(X;\H_0)=C_q^*(X;\H_0)$. So we only focus on the converse.

Assume that $X$ does not have Property~A. We aim to construct an operator $T \in C^*_{uq}(X) \setminus C^*_u(X)$, which will conclude the proof of Theorem~\ref{introthm:uniform-converse} by Corollary~\ref{cor:scalar-reductions}(1). In fact, our proof below shows that $T\in C^*_{uq}(X) \setminus C^*_u(X)$ can moreover be chosen to be a ghost. 
Let us fix a sequence of subspaces $\{Y_n\}_{n\in \NN}$ satisfying \eqref{EQ:sparse} and a sequence of positive contractions $\{S_n\}_n$ with $S_n \in \B(\ell^2(Y_n))$ by Proposition~\ref{prop:weak expander spectral}, satisfying the conditions therein. Set $q(s):=\sup_n q_s(S_n)$ for $s>0$. 

For each integer $k\geq1$, set $d_k:=2^k$ and choose $r_k\geq k$
such that $q(r_k)^2\leq\frac{1}{kd_k}$. Set
\[
 \gamma_k:=\frac{1}{kd_k^2}\quad \text{and} \quad
 \tau_k:=\frac{1}{kN(r_k+1)}.
\]
Let $\ell_k>0$ be as in Proposition~\ref{prop:range-split} for
$k,\tau_k,\gamma_k$. By condition (3) in Proposition~\ref{prop:weak expander spectral}, we may choose
strictly increasing $n(k)\in\NN$ such that 
\[
\max_{x\in Y_{n(k)}}\|S_{n(k)}\delta_x\|^2\leq \min\left\{\ell_k, \frac{1}{k^2}\right\}.
\] 
Simply denote $Z_k:=Y_{n(k)}$ and $R_k:=S_{n(k)}$.
Then $R_k$ is a positive contraction with
\[
 \sum_{x\in Z_k}\theta_{R_k\delta_x}=R_k\left(\sum_{x\in Z_k}\theta_{\delta_x}\right)R_k^*
 =R_k^2\leq\Id_{\ell^2(Z_k)}.
\]
Apply Proposition~\ref{prop:range-split} to $Y=Z_k$, $\H=\ell^2(Z_k)$, $u_x:=R_k\delta_x$ for $x\in Z_k$ and the parameters $\gamma_k,\tau_k$. This provides pairwise disjoint $C_\sigma^{(k)}\subseteq Z_k$ for $\sigma\in\{-1,1\}^k$, satisfying the conditions therein. Put
\[
 G_k:=\bigcup_\sigma C_\sigma^{(k)} 
 \quad \text{and} \quad B_k:=Z_k\setminus G_k.
\]
Then $|B_k|\leq\tau_k|Z_k|$ by Proposition~\ref{prop:range-split}. Define
\[
 \Lambda_k:=\N_{r_k+1}(B_k) \cap Z_k,
 \quad Q_k:=\chi_{Z_k\setminus\Lambda_k}
 \quad \text{and} \quad H_k:=R_kQ_k.
\]
Then $H_k$ is a contraction, $\|\chi_{B_k}H_k\|\leq q(r_k)$ and $|\Lambda_k|\leq N(r_k+1)|B_k|\leq N(r_k+1)\tau_k|Z_k|$.

For any self-adjoint $A\in\B(\ell^2(Z_k))$, applying
\eqref{eq:range-cells} to $Q_kAQ_k$, we obtain
\[
 \sum_{x\in C_\sigma^{(k)}}|\langle Q_kAQ_kR_k\delta_x,R_k \delta_x\rangle|
 \leq d_k^{-1}\sum_{x\in Z_k}|\langle Q_kAQ_kR_k\delta_x,R_k \delta_x\rangle|
       +\gamma_k\|Q_kAQ_k\|_1.
\]
Since $Q_kR_k\delta_x = H_k^*\delta_x$ and $\|Q_kAQ_k\|_1\leq\|A\|_1$, it follows that the assumptions of Proposition~\ref{prop:constructing quasi-local oeprators} hold for $Y=Z_k$, $H=H_k$, $C_\sigma=C_\sigma^{(k)}$ and
$\gamma=\gamma_k$. Moreover,
\[
 \varepsilon_k:=d_k\|\chi_{B_k}H_k\|^2+d_k^2\gamma_k
 \leq d_kq(r_k)^2+d_k^2\gamma_k\leq\frac2k.
\]
By the construction from Proposition~\ref{prop:constructing quasi-local oeprators}, we obtain operators $V^{(k)}_j \in \B(\ell^2(Z_k))$ and $V^{(k)}:\CC^k\otimes\ell^2(Z_k)\to\ell^2(Z_k)$. We record their useful properties as follows.

\begin{lem}\label{lem:constructing quasi-local candidates}
The operators $V^{(k)}$ satisfy the following: 
\begin{enumerate}
 \item $M:=\sup_k\|V^{(k)}\|<\infty$;
 \item $\omega(s):= \sup_{k\geq1}\sup_{\substack{E,D\subseteq Z_k\\d(E,D)>s}}
\|\chi_EV^{(k)}(\Id_{\CC^k}\otimes\chi_D)\|\to 0$ as $s \to +\infty$;\\[.2cm]
 \item $b_k:=\frac1k\left\|\sum_{j=1}^k(V_j^{(k)})^*V_j^{(k)}\right\|\to 1$ as $k\to+\infty$;
 \item $\|(V^{(k)})^*\delta_x\|^2 \leq \frac1k$ for each $k \in \NN$ and $x\in Z_k$.
\end{enumerate}
\end{lem}

\begin{proof}
Given any $E,D \subseteq Z_k$, note that
\[
\|\chi_EH_k\chi_D\| = \|\chi_ER_kQ_k\chi_D\| = \|\chi_ER_k\chi_DQ_k\| \leq \|\chi_ER_k\chi_D\|.
\]
Hence, \eqref{eq:abstract-row-bound} in Proposition~\ref{prop:constructing quasi-local oeprators} shows that
\[
 \|\chi_EV^{(k)}(\Id_{\CC^k}\otimes\chi_D)\|^2
 \leq h^2+\sqrt2\,h+\frac4k+\frac2{\sqrt{k}}
 \quad \text{for} \quad h:=\|\chi_ER_k\chi_D\|.
\]
In particular, condition (1) holds. To see (2), note that $h=\|\chi_ER_k\chi_D\|\leq q(s)$ for any $d(E,D)>s$ by definition. Given $\delta>0$, choose $K$ such that $4/k+2/\sqrt{k}<\delta^2/2$ for any $k\geq K$. Then choose $s>0$ so that
$q(s)^2+\sqrt2q(s)<\delta^2/2$ and $s>\max_{1\leq k<K}\diam(Z_k)$ (when $K=1$, the latter
condition is omitted). Thus $\omega(s)\leq\delta$, which proves (2).

For (3), set $e_k:=|Z_k|^{-1/2}\chi_{Z_k}$. Since $R_ke_k=e_k$, we have
\[
 \|(\Id-Q_k)e_k\|^2
 =\frac{|\Lambda_k|}{|Z_k|}
 \leq N(r_k+1)\tau_k\leq\frac1k
 \quad \text{and} \quad 
 \|\chi_{B_k}e_k\|^2\leq\tau_k.
\]
Hence, $\|\chi_{G_k}H_ke_k-e_k\|
 \leq\|(\Id-Q_k)e_k\|+\|\chi_{B_k}e_k\|
 \leq k^{-1/2}+\sqrt{\tau_k}\to 0$. 
By construction, we have $(V_j^{(k)})^*V_j^{(k)}=H_k^*\chi_{G_k}H_k$ for each $j$. Therefore, $b_k =\|\chi_{G_k}H_k\|^2\to 1$ as $k \to +\infty$. 

Finally for (4), note that
\[
\|(V^{(k)})^*\delta_x\|^2=\sum_{j=1}^k\|(V_j^{(k)})^*\delta_x\|^2\le k\|H_k^*\delta_x\|^2=k\|Q_kR_k\delta_x\|^2\le k\|R_k\delta_x\|^2\le\frac1k.
\]
So we conclude the proof.
\end{proof}

Given a unitary $U$ on $\CC^k$, set
\[
 T_U^{(k)}:=V^{(k)}(U\otimes\Id_{\ell^2(Z_k)})(V^{(k)})^* \in \B(\ell^2(Z_k)).
\]
Lemma~\ref{lem:constructing quasi-local candidates}(1) shows that $\|T_U^{(k)}\|\leq M^2$, and it is easy to deduce from Lemma~\ref{lem:constructing quasi-local candidates}(2) that $q_{3s}(T_U^{(k)})\leq2M\omega(s)$ for any $s>0$. For any \(x,y\in Z_k\), Lemma~\ref{lem:constructing quasi-local candidates}(4) shows
\begin{equation}\label{EQ:ghost estimate}
\left|(T_U^{(k)})_{x,y}\right|=|\langle (U\otimes \Id_{\ell^2(Z_k)})(V^{(k)})^*\delta_y,
(V^{(k)})^*\delta_x
\rangle|\le
\|(V^{(k)})^*\delta_y\| \cdot 
\|(V^{(k)})^*\delta_x\|\le\frac1k.
\end{equation}
Define
\[
 F_k:\ell^2(Z_k)\longrightarrow\ell^2(Z_k;\CC^k)\quad \text{by} \quad
 F_k\xi:=\frac1{\sqrt k}
 \bigl(V_1^{(k)}\xi,\ldots,V_k^{(k)}\xi\bigr).
\]
Then Lemma~\ref{lem:constructing quasi-local candidates}(3) shows that $\|F_k\|^2 = b_k \to 1$ as $k\to+\infty$.

Take $\Gamma_k$ to be the symmetric group on $k+1$ letters, and $\pi_k:\Gamma_k \to \B(\CC^k)$ to be its standard irreducible unitary representation. For each $g\in \Gamma_k$, denote by $\overline{\pi_k(g)}$ the entry-wise complex conjugate of $\pi_k(g)$. We have the following.

\begin{lem}\label{lem:Schur orthogonality}
We have $F_kF_k^* = \frac1{|\Gamma_k|}\sum_{g\in\Gamma_k} T_{\pi_k(g)}^{(k)}\otimes\overline{\pi_k(g)}$.
\end{lem}

\begin{proof}
By Schur orthogonality, we have 
\[
\frac1{|\Gamma_k|}\sum_{g\in\Gamma_k}
\pi_k(g)_{i,j}\overline{\pi_k(g)_{p,q}}
=\frac1k\delta_{ip}\delta_{jq} \quad \text{for any} \quad 1 \leq i,j,p,q \leq k,
\]
where $\delta_{ij}$ denotes the Kronecker delta. Hence for $1 \leq p,q \leq k$, we have
\begin{align*}
\left(
\frac1{|\Gamma_k|}\sum_{g\in\Gamma_k}
T_{\pi_k(g)}^{(k)}\otimes\overline{\pi_k(g)}
\right)_{p,q}
&=\frac1{|\Gamma_k|}\sum_{g\in\Gamma_k}
T_{\pi_k(g)}^{(k)}\overline{\pi_k(g)_{p,q}}\\
&=\sum_{i,j=1}^k V^{(k)}_i(V^{(k)}_j)^*
\left(
\frac1{|\Gamma_k|}\sum_{g\in\Gamma_k}
\pi_k(g)_{i,j}\overline{\pi_k(g)_{p,q}}
\right)\\
&=
\frac1k\sum_{i,j=1}^k
\delta_{ip}\delta_{jq}V^{(k)}_i(V^{(k)}_j)^*\\
&=\frac1kV_p^{(k)}(V^{(k)}_q)^*=(F_kF_k^*)_{p,q}.
\end{align*}
This concludes the proof.
\end{proof}


Now we adapt Ozawa's argument in \cite[Lemma 1]{Oza} to construct an operator in $C^*_{uq}(X) \setminus C^*_u(X)$. Given $\rho>0$, denote $\MM_\rho(Z_k):=\{T \in \B(\ell^2(Z_k)): \ppg(T) \leq \rho\}$. 

\begin{lem}\label{lem:Ozawa}
Given $\rho>0$, there exists $K(\rho)\in\NN$ such that for every $k\geq K(\rho)$ there exists $g\in \Gamma_k$ satisfying
\begin{equation}\label{EQ:distance estimate}
 \dist\bigl(T_{\pi_k(g)}^{(k)},\MM_\rho(Z_k)\bigr)\geq\tfrac14.
\end{equation}
\end{lem}

\begin{proof}
Fix $\rho>0$, and assume the result fails for $\rho$. Then there are infinitely many $k\in \NN$ such that for any $g\in \Gamma_k$ we can find $c_g\in \MM_\rho(Z_k)$ satisfying $\big\|T_{\pi_k(g)}^{(k)}-c_g\big\|<1/4$. Then $\|c_g\|\leq M^2+\frac{1}{4}$, and Lemma~\ref{lem:Schur orthogonality} gives
\begin{equation}\label{eq:lower}
 \left\|\frac1{|\Gamma_k|}\sum_{g\in\Gamma_k}
       c_g\otimes\overline{\pi_k(g)}\right\| \geq b_k-\tfrac14.
\end{equation}
Recall that $\{(x,y)\in Z_k^2:d(x,y)\leq \rho\}$ can be written as a union of at most $2N(\rho)$ graphs of partial translations (see, \emph{e.g.}, \cite[Lemma~2]{Oza}). Consequently, whenever $b_k>1/4$,
\eqref{eq:lower} yields some $x,y\in Z_k$ such that
\[
\left\|\frac1{|\Gamma_k|}\sum_{g\in\Gamma_k}
  (c_g)_{x,y}\,\overline{\pi_k(g)}\right\| \geq \frac{b_k-1/4}{2N(\rho)}.
\]
On the other hand, \cite[Lemma~3]{Oza} shows that for $x,y$ above, we have
\[
\left\|\frac1{|\Gamma_k|}\sum_{g\in\Gamma_k}
  (c_g)_{x,y}\,\overline{\pi_k(g)}\right\| \leq \frac{\sup_{g\in \Gamma_k}\left|(c_g)_{x,y}\right|}{\sqrt k} \leq \frac{M^2+1/4}{\sqrt k}.
\]
This is impossible for sufficiently large $k$, since $b_k \to 1$ as $k \to +\infty$. 
\end{proof}

Finally, choose strictly increasing positive integers $\{k_m\}_{m \in \NN}$ with
$k_m\geq K(m)$ from Lemma~\ref{lem:Ozawa}, and for each $m\geq1$ choose
$g_m\in \Gamma_{k_m}$ satisfying \eqref{EQ:distance estimate} for $\rho=m$. Define
\[
T:=\mathrm{(SOT)}\text{-}\sum_{m\in \NN} T_{\pi_{k_m}(g_m)}^{(k_m)} \in \prod_m \B(\ell^2(Z_{k_m})) \subseteq \B(\ell^2(X)).
\]
Then $q_{3s}(T)\leq2M\omega(s)\to 0$ as $s\to+\infty$, which means that $T \in C^*_{uq}(X)$. For any finite-propagation operator $S'$ on $\ell^2(X)$, choose $m\geq1$ with $m\geq\ppg(S')$. Then $\chi_{Z_{k_m}}S'\chi_{Z_{k_m}}$ belongs to $\MM_m(Z_{k_m})$, whence $\|T-S'\|\geq1/4$. Hence $\dist(T,C_u^*(X))\geq1/4$, and so $T \in C^*_{uq}(X) \setminus C^*_u(X)$. 

This concludes the proof of Theorem~\ref{introthm:uniform-converse}.
Moreover, \eqref{EQ:ghost estimate} shows that
\[
\sup_{x,y\in Z_{k_m}}|T_{x,y}|
\le\frac1{k_m}\longrightarrow0 \quad \text{as} \quad m \to +\infty.
\]
Hence, we obtain the following stronger version of the converse in Theorem~\ref{introthm:uniform-converse}.

\begin{prop}\label{prop:ghost}
Let $X$ be a discrete metric space of bounded geometry, which does not have Property~A. Then there exists a ghost operator $T \in C^*_{uq}(X) \setminus C^*_u(X)$.
\end{prop}

\section{Maximal uniform Roe algebras}\label{sec:maximal}

In this section, we prove Theorem~\ref{introthm:maximal-converse},
which characterises weak containment for Roe algebras
in terms of Property~A. As in the proof of Theorem~\ref{introthm:uniform-converse}, it suffices to consider the scalar case $\H_0=\CC$ by Corollary~\ref{cor:scalar-reductions}(2). We will use the operators from the failure of Property~A again (Proposition~\ref{prop:weak expander spectral}) together with the signing estimate from Lemma~\ref{lem:range-signing} and the tools of random graphs.

\subsection{Random graphs and signed adjacency matrices}\label{subsec:random-graphs}

Let us recall the graph-theoretic notation used in our construction. Given a simple (\emph{i.e.}, without loops or multiple edges) graph $G=(V,E)$ and a vertex $x\in V$, its \emph{degree} $\deg_G(x)$ is the number of edges with one endpoint $x$. The \emph{mean degree} of $G$ is 
\[
\mdeg(G):=\frac{1}{|V|} \sum_{x\in V} \deg_G(x) = \frac{2|E|}{|V|}.
\]
The \emph{girth} of $G$, denoted by $\girth(G)$, is the length of a shortest cycle in $G$. When $G$ contains no cycles, set $\girth(G):=+\infty$. Also recall that its \emph{adjacency matrix} \(A_G\in\B(\ell^2(V))\) is defined by
\[
(A_G)_{x,y}:=
\begin{cases}
1,&\{x,y\}\in E,\\
0,&\text{otherwise}.
\end{cases}
\]


The following result produces the graphs we need from a non-Property~A space.

\begin{prop}\label{prop:large-girth-graphs}
Let $X$ be a discrete metric space of bounded geometry, which does not have Property~A. For each integer $b\ge2$, there
exist pairwise disjoint non-empty finite subsets $Y_n\subseteq X$ and
finite simple graphs $G_n=(Y_n, E_n)$ satisfying the following:
\begin{enumerate}
\item $\sup_{x\in Y_n} \deg_{G_n}(x) \leq 2b$, and $\mdeg(G_n) > b-2$;
\item $\girth(G_n)\to\infty$ as $n \to +\infty$;
\item there exists $R>0$ which only depends on $b$ and $X$, such that
      $d(x,y)\le R$ for all $n$ and $\{x,y\}\in E_n$.
\end{enumerate}
\end{prop}

We prove Proposition~\ref{prop:large-girth-graphs} by selecting edges
at random and then deleting edges to remove short cycles and excess
degrees. We first recall some elementary facts about random graphs. Let $V$ be a non-empty finite set, and put
\[
 \mathcal E(V)=\big\{\{x,y\}:x,y\in V,\ x\ne y\big\}
 \quad \text{and} \quad \Omega=2^{\mathcal E(V)}.
\]
For each $e\in\mathcal E(V)$, fix $p_e\in[0,1]$.  Define a probability measure
$\mathbb P$ on $2^\Omega$ by
\[
 \mathbb P(\{\omega\})=
 \prod_{e\in\omega}p_e
 \prod_{e\in\mathcal E(V)\setminus\omega}(1-p_e)
 \quad \text{for} \quad \omega\in\Omega.
\]
Each $\omega\in\Omega$ gives a
simple graph $F_\omega=(V,\omega)$.  We call $F=\{F_\omega: \omega \in \Omega\}$ a
\emph{random simple graph} on $V$ with independent edge probabilities
$\{p_e:e\in \mathcal E(V)\}$.  For a function $f:\Omega\to\mathbb R$, denote the expectation and variance by
\[
 \mathbb Ef:=\int_\Omega f\,\mathrm{d}\mathbb P
 \quad \quad \text{and} \quad \operatorname{Var}(f):=\mathbb E[(f-\mathbb Ef)^2].
\]
For $e\in\mathcal E(V)$, denote
$A_e:=\{\omega\in\Omega:e\in\omega\}$.  For any distinct $e_1,\ldots,e_k\in\mathcal E(V)$, we have
\begin{equation}\label{eq:edge-products}
 \mathbb E\left(\prod_{i=1}^k \chi_{A_{e_i}}\right) = \prod_{i=1}^k \mathbb E\left(\chi_{A_{e_i}}\right) = \prod_{i=1}^k \mathbb P(A_{e_i})=\prod_{i=1}^k p_{e_i}.
\end{equation}

\begin{lem}\label{lem:random-graph-calculations}
Let $F$ be a random simple graph on $V$ with edge probabilities $\{p_e:e\in \mathcal E(V)\}$. Let $Q=(q_{x,y})_{x,y\in V}$ be given by $q_{x,y}=p_{\{x,y\}}$ for
$x\ne y$ and $q_{x,x}=0$. For $x\in V$, set
$D_x(\omega)=\deg_{F_\omega}(x)$. For an integer $j\ge3$, 
let $C_j(\omega)$ be the number of simple cycles of length $j$ in
$F_\omega$. Then
\[
 \mathbb ED_x=\sum_{y\in V}q_{x,y},\quad
 \operatorname{Var}(D_x)=\sum_{y\in V}q_{x,y}(1-q_{x,y})
       \le\mathbb ED_x
 \quad \text{and} \quad 
 \mathbb EC_j\le\frac{\Tr(Q^j)}{2j}.
\]
\end{lem}

\begin{proof}
Fixing $x\in V$, we have $D_x=\sum_{y\ne x}\chi_{A_{\{x,y\}}}$. Then \eqref{eq:edge-products} gives 
\[
\mathbb ED_x=\sum_{y \neq x}\mathbb E\chi_{A_{\{x,y\}}} =\sum_{y\in V}q_{x,y}
\] 
and the independence of the edges gives
\begin{align*}
\operatorname{Var}(D_x) &=\sum_{y\ne x}\operatorname{Var}(\chi_{A_{\{x,y\}}})=\sum_{y\ne x}
\left(\mathbb E\chi_{A_{\{x,y\}}}^2
-(\mathbb E\chi_{A_{\{x,y\}}})^2\right)=\sum_{y\ne x}(p_{\{x,y\}}-p_{\{x,y\}}^2)\\
&=\sum_{y\in V}q_{x,y}(1-q_{x,y})\le\sum_{y\in V}q_{x,y}
=\mathbb ED_x.
\end{align*}
For the cycles, each simple cycle of length $j$ is counted $2j$
times by ordered choices of an initial vertex and an orientation.
Thus, writing $x_j=x_0$,
\[
 \mathbb EC_j
 =\frac1{2j}
   \sum_{\substack{x_0,\ldots,x_{j-1}\in V\\
                    \text{pairwise distinct}}}
     \prod_{i=0}^{j-1}q_{x_i,x_{i+1}}\le\frac1{2j}
   \sum_{x_0,\ldots,x_{j-1}\in V}
     \prod_{i=0}^{j-1}q_{x_i,x_{i+1}}
 =\frac{\Tr(Q^j)}{2j}.
\]
Here the inequality uses the non-negativity of the matrix entries.
\end{proof}

Based on the calculations in Lemma~\ref{lem:random-graph-calculations}, we have the following.

\begin{lem}\label{lem:graph-selection}
Let $V$ be a non-empty finite set and $Q=(q_{x,y})_{x,y\in V}$ be a symmetric matrix with $q_{x,x}=0$ and $0\le q_{x,y}\le1$.
Suppose integers $b\ge2$ and $c\ge3$ satisfy
\[
 b-\tfrac12\le\sum_{y\in V}q_{x,y}\le b
 \quad \text{for each} \quad x\in V, \quad \text{and} \quad 
 \sum_{j=3}^c\frac{\Tr(Q^j)}{2j}\le\frac{|V|}{4}.
\]
Then there is a simple graph $G$ on $V$ with $\sup_{x\in V} \deg_{G}(x) \leq 2b$, $\mdeg(G) > b-2$, $\girth(G)>c$ such that $E(G)\subseteq\{\{x,y\}:x\ne y,\ q_{x,y}>0\}$.
\end{lem}

\begin{proof}
Define the random simple graph $F$ on $V$ with edge probabilities $q_{x,y}$. Let $D_x$ and $C_j$ be as in
Lemma~\ref{lem:random-graph-calculations}. Then Lemma~\ref{lem:random-graph-calculations} yields
\[
 b-\tfrac12\le \mathbb ED_x\le b,\quad
 \operatorname{Var}(D_x)\le \mathbb ED_x
 \quad \text{and} \quad \mathbb E\left(\sum_{j=3}^c C_j\right)\le\frac{|V|}{4}.
\]
For $u\in\mathbb R$, write $u_+=\max\{u,0\}$.
Applying the inequality $(u-a)_+\le \frac{u^2}{4a}$ for $a>0$ to 
$u=D_x-\mathbb ED_x$ and $a=2b-\mathbb ED_x\ge b$, we obtain
\[
 \mathbb E(D_x-2b)_+
 \le\frac{\operatorname{Var}(D_x)}{4(2b-\mathbb ED_x)}
 \le\frac14.
\]
Denote the function $|E(F)(\omega)|:=|E(F_\omega)|$ for $\omega \in \Omega$. Then
\[
 \frac{2}{|V|}\mathbb E\left(
 |E(F)|-\sum_{j=3}^{c}C_j
           -\sum_{x\in V}(D_x-2b)_+\right)\ge\frac{1}{|V|}\sum_{x\in V}\mathbb ED_x
            -\frac{2}{|V|}\left(\frac{|V|}{4}+\frac{|V|}{4}\right)
 \ge b-\frac32>b-2.
\]

Hence, there is $\omega\in\Omega$ with $\mathbb P(\{\omega\})>0$ for which 
\begin{equation}\label{EQ:graph deletion}
|E(F_\omega)|-\sum_{j=3}^{c}C_j(\omega) -\sum_{x\in V}(D_x(\omega)-2b)_+ > \frac{|V|(b-2)}{2}.
\end{equation}
Repeatedly delete an edge from a cycle of length at most \(c\) until no such cycle remains. This deletes at
most $\sum_{j=3}^c C_j(\omega)$ edges.  Next, as long as some vertex has degree greater than $2b$, delete an edge incident to it. Each such deletion decreases the sum of excess degrees by at least one, so at
most $\sum_{x\in V}(D_x(\omega)-2b)_+$ further edges are removed.
The resulting graph $G$ has the required degree and girth bounds,
and \eqref{EQ:graph deletion} gives its mean-degree bound.
Finally, $\mathbb P(\{\omega\})>0$ ensures that every edge of
$F_\omega$, hence of $G$, has $q_{x,y}>0$.
\end{proof}

We apply the graph-selection lemma to powers of the random-walk
operators arising from the local spectral estimates in
Section~\ref{sec:preparations}.

\begin{proof}[Proof of Proposition~\ref{prop:large-girth-graphs}]
Fix an integer $b \geq 2$. We use the construction from the proof of Proposition~\ref{prop:weak expander spectral}. More precisely, fix $r_0,h>0$, a sequence of subspaces $\{Y_n\}_{n\in \NN}$ from Proposition~\ref{prop:pieces} and set $D:=N(r_0)$. For each $n\in \NN$, let $H_n:=G_{r_0}(Y_n)$, $\Delta_n$ be the graph Laplacian, and put $L_n:=\frac{\Delta_n}{2D}$ and $P_n:=\Id_{\ell^2(Y_n)}-L_n$. Also set $\rho:=1-\frac{h^2}{4D^2}$ and $t_n=\left\lfloor\frac{n}{2r_0}\right\rfloor$. Then $0\le P_n\le \Id_{\ell^2(Y_n)}$ and clearly, $P_n$ is symmetric and stochastic (\emph{i.e.}, its entries are non-negative and each row sums to one). Applying Proposition~\ref{prop:pieces} to a singleton gives $0<h\le D$ and hence, $0<\rho<1$. 

Claim: for $s=1,\cdots,t_n$ and $x\in Y_n$, we have $(P_n^s)_{x,x}\le\rho^s$. Indeed, for any non-empty $B\subseteq Y_n$ with $\diam(B)\le n$, the second paragraph in the proof of Proposition~\ref{prop:weak expander spectral} shows that the spectrum of $\chi_B L_n \chi_B \in \B(\ell^2(B))$ lies in $[1-\rho,1]$ and hence, $0\le \chi_B P_n \chi_B\le\rho \Id_{\ell^2(B)}$. Now fix $x\in Y_n$ and take $B:=\{y\in Y_n:d_{H_n}(x,y)\le t_n\}$, where $d_{H_n}$ is the graph distance. Then $\diam_X(B)\le2r_0t_n \le n$ and for each integer $1\le s\le t_n$, we have $(P_n^s)_{x,x}=((\chi_BP_n\chi_B)^s)_{x,x} \le\|\chi_BP_n\chi_B\|^s\le\rho^s$.

Now choose an even integer $k\ge2$ such that $b\rho^k\le\frac12$.
Discarding finitely many indices and relabelling, we assume that $t_n\ge3k$ for all $n$. We have $c_n:=\lfloor t_n/k\rfloor \to \infty$. For each $n$, define a symmetric matrix
$Q_n:=(q_{x,y}^{(n)})_{x,y\in Y_n}$ by
\[
 q_{x,y}^{(n)}=
 \begin{cases}
 b(P_n^k)_{x,y},&x\ne y;\\
 0,&x=y.
 \end{cases}
\]
Since $k$ is even, the Cauchy--Schwarz inequality and
the claim above give
\[
0\le (P_n^k)_{x,y}=\big|\langle P_n^{k/2}\delta_y,
                 P_n^{k/2}\delta_x\rangle\big|\le \|P_n^{k/2}\delta_y\| \cdot 
      \|P_n^{k/2}\delta_x\|=\big((P_n^k)_{x,x}(P_n^k)_{y,y}\big)^{1/2}
\le \rho^k
\]
for $x,y\in Y_n$, which shows that $0\le q_{x,y}^{(n)}\le\frac12$.  As $P_n^k$ is stochastic, we have
\[
 b-\tfrac12\le
 \sum_{y\in Y_n}q_{x,y}^{(n)}
 =b\bigl(1-(P_n^k)_{x,x}\bigr)\le b.
\]
For each integer $3\le j\le c_n$, we have $kj\le t_n$. From the claim above, we obtain
\[
 \Tr(Q_n^j)\le b^j\Tr(P_n^{kj})
 \le |Y_n|(b\rho^k)^j.
\]
Hence
\[
 \sum_{j=3}^{c_n}\frac{\Tr(Q_n^j)}{2j}
 \le |Y_n|\sum_{j=3}^\infty 2^{-j}
 =\frac{|Y_n|}{4}.
\]
Applying Lemma~\ref{lem:graph-selection} with $V=Y_n$ and $c=c_n$
yields graphs $G_n$ satisfying (1) and (2). If
$\{x,y\}\in E(G_n)$, then $q_{x,y}^{(n)}>0$ and hence,
$(P_n^k)_{x,y}>0$. By the definition of $P_n$, this implies that
$x$ and $y$ are joined by an $H_n$-path of length at most $k$, so
$d(x,y)\le kr_0$. Thus (3) holds with $R=kr_0$.
\end{proof}

For the unsigned adjacency matrix $A_{G_n}$ of $G_n$ constructed in Proposition~\ref{prop:large-girth-graphs}. With $\xi_n:=|Y_n|^{-1/2}\chi_{Y_n}$, we have
\[
b-2<\frac{2|E(G_n)|}{|Y_n|}
=\langle A_{G_n}\xi_n,\xi_n\rangle\le\|A_{G_n}\|
\le\max_{x\in Y_n}\deg_{G_n}(x)
\le2b.
\]
To reduce the norm of the adjacency matrix, we consider an edge signing. More precisely, given an edge signing \(\varepsilon:E\to\{-1,1\}\), the corresponding \emph{signed adjacency matrix} \(A_{G,\varepsilon}\) is defined by
\[
(A_{G,\varepsilon})_{x,y}:=
\begin{cases}
\varepsilon(\{x,y\}),&\{x,y\}\in E,\\
0,&\text{otherwise}.
\end{cases}
\]
The signing estimate in Lemma~\ref{lem:range-signing} now allows
us to reduce the norm of the adjacency matrix from order $b$ to order
$\sqrt b$.

\begin{lem}\label{lem:signed-adjacency}
Let $G=(V,E)$ be a finite simple graph with maximum degree at most
$\Delta$, where $\Delta\ge1$.  There exists an edge signing
$\varepsilon:E\to\{-1,1\}$ such that its signed adjacency matrix
$A_{G,\varepsilon}\in\B(\ell^2(V))$ satisfies $\|A_{G,\varepsilon}\|\le32\sqrt\Delta$.
\end{lem}

\begin{proof}
The result holds trivially when $E = \emptyset$. Otherwise, for each $e=\{x,y\}\in E$, set 
\[
 u_e:=\frac{\delta_x+\delta_y}{\sqrt{2\Delta}}
 \in\ell^2(V)
\]
and consider the rank-one operator $\theta_{u_e} \in \B(\ell^2(V))$ from \eqref{EQ:theta operator}. Since the degree of each vertex is at most $\Delta$, we have 
\begin{align*}
\sum_{e\in E}|\langle f,u_e\rangle|^2&=\frac1{2\Delta}
\sum_{\{x,y\}\in E}|f(x)+f(y)|^2\le\frac1{\Delta}
\sum_{\{x,y\}\in E}
\bigl(|f(x)|^2+|f(y)|^2\bigr)\\
&=\frac1{\Delta}
\sum_{x\in V}\deg_G(x)|f(x)|^2\le\sum_{x\in V}|f(x)|^2
=\|f\|^2 \quad \text{for} \quad f\in\ell^2(V).
\end{align*}
Hence, $\sum_e\theta_{u_e}\le \Id_{\ell^2(V)}$ and $\|u_e\|^2=1/\Delta$.
Applying Lemma~\ref{lem:range-signing} with $d=1$ and
$\ell=1/\Delta$, we obtain signs $\varepsilon_e\in\{-1,1\}$ such that
\[
 \|B\|\le\frac8{\sqrt\Delta} \quad \text{for} \quad B:=\sum_{e\in E}\varepsilon_e\theta_{u_e}.
\]
Let $\operatorname{diag}(B)$ be the diagonal part of $B$ relative to the basis $(\delta_x)_{x\in V}$. Direct calculations show that $A_{G,\varepsilon}=2\Delta(B-\operatorname{diag}(B))$.  Hence $\|A_{G,\varepsilon}\|\le4\Delta\|B\|\le32\sqrt\Delta$.
\end{proof}

\subsection{Signed distance kernels and maximal norm estimates}\label{subsec:maximal-norm}

We next prove a lower bound for the maximal norm of signed
adjacency matrices of graphs with large girth.

\begin{prop}\label{prop:maximal-norm-lower}
Let $X$ be a discrete metric space of bounded geometry, and $G_n=(Y_n, E_n)$ be finite simple graphs on pairwise disjoint non-empty subsets $Y_n\subseteq X$. Suppose that the degrees of $G_n$ are uniformly
bounded, $\girth(G_n)\to\infty$, and there is $R>0$ such that $d(x,y)\le R$ for every $n$ and $\{x,y\}\in E_n$.
For each $n$, fix an edge signing $\varepsilon_n:E_n\to\{-1,1\}$ and let
$A_n\in\B(\ell^2(Y_n))$ be the corresponding signed adjacency matrix. Take $A:=\mathrm{(SOT)}$-$\sum_n A_n \in \B(\ell^2(X))$. Then $A=A^*\in\CC_u[X]$ and
\[
 \|A\|_{u,\max}\ge
 \limsup_{n\to\infty}\frac{2|E_n|}{|Y_n|}.
\]
\end{prop}

We need signed distance kernels for the proof. For the graphs $G_n$
above, let $d_n$ be the graph distance on $Y_n$, which takes $+\infty$ between different connected components. For each unordered pair $\{x,y\}\subseteq Y_n$ in the same component, choose a shortest path $P_{xy}$ from $x$ to $y$, with $P_{yx}$ its reversal. For a path $P$, denote by $s(P)$ the product of the signs of its edges, with $s(P)=1$ for the empty path. Given $0<q<1$, define a real symmetric matrix $K_n^{(q)}\in\B(\ell^2(Y_n))$ by
\[
 K_n^{(q)}(x,y)=
 \begin{cases}
 q^{d_n(x,y)}s(P_{xy}),& \text{if }d_n(x,y)<\infty,\\
 0,&\text{if }d_n(x,y)=\infty.
 \end{cases}
\]
In particular, $K_n^{(q)}(x,x)=1$ and $K_n^{(q)}(x,y)=q\varepsilon_n(\{x,y\})$ when
$\{x,y\}\in E_n$.

\begin{lem}\label{lem:signed-kernel}
There exist positive integers $a_n\to\infty$ such that for any integer \(m\ge1\), there is \(N_m\in\mathbb N\) satisfying
\[
 K_n^{(q)}|_{Z\times Z}\ge-2mq^{a_n}\Id_{\ell^2(Z)}
\]
for any $0<q<1$, $n \geq N_m$ and any $Z\subseteq Y_n$ with $|Z|\le m$.
\end{lem}

\begin{proof}
Set \(c_n:=\operatorname{girth}(G_n)\) and $a_n:=\left\lfloor\sqrt{\min\{n,c_n\}}\right\rfloor$. Fix an integer \(m\ge1\). Since \(c_n\to\infty\) and \(a_n\le\sqrt{c_n}\), there exists \(N_m\in\mathbb N\) such that $m^2a_n<c_n$ for $n\ge N_m$. Fix \(n\ge N_m\) and a non-empty subset \(Z\subseteq Y_n\) with \(|Z|\le m\). Let \(F\) be the subgraph of \(G_n\) obtained by taking the union of the chosen shortest paths \(P_{xy}\) for all \(x,y\in Z\) with \(d_n(x,y)\le a_n\), together with the remaining vertices of \(Z\). Since $|E(F)|\le m^2a_n<c_n$, the graph \(F\) contains no cycles. Hence every connected component of \(F\) is a tree.

For each connected component \(T\) of \(F\), choose a root \(o\in V(T)\). Denote by \([x,y]\) the unique path in \(T\) joining \(x,y\in V(T)\), and set $\sigma_x:=s([o,x])\in\{-1,1\}$ for $x\in V(T)$. Let \(d_F\) denote the graph distance in \(F\). Since tree distances are conditionally negative definite, Schoenberg's theorem (\cite[Theorem~1]{Sch}) implies that $(q^{d_F(x,y)})_{x,y\in V(T)}$ is positive semidefinite for every \(0<q<1\). Moreover, $\sigma_x\sigma_y=s([x,y])$ for every $x,y\in V(T)$. Conjugating by the diagonal matrix with entries \(\sigma_x\), we see that $(\sigma_x\sigma_yq^{d_F(x,y)})_{x,y\in V(T)}$ is also positive semi-definite.

For each \(0<q<1\), define \(L^{(q)}\in\B(\ell^2(Z))\) by
\[
L^{(q)}_{x,y}:=
\begin{cases}
\sigma_x\sigma_yq^{d_F(x,y)},
&\text{if } x,y \text{ belong to the same component of } F;\\
0,&\text{otherwise}.
\end{cases}
\]
By the preceding argument, we have \(L^{(q)}\ge0\). We now compare \(L^{(q)}\) with \(K_n^{(q)}|_{Z\times Z}\). 

If \(x,y\in Z\) satisfy \(d_n(x,y)\le a_n\), then \(P_{xy}\) is contained in \(F\). Since \(F\) contains no cycles, \(P_{xy}\) is the unique path joining \(x\) and \(y\) in \(F\). Hence, $d_F(x,y)=d_n(x,y)$ and $\sigma_x\sigma_y=s(P_{xy})$, which shows that $L^{(q)}(x,y)=K_n^{(q)}(x,y)$. 
If \(d_n(x,y)>a_n\) (including \(d_n(x,y)=+\infty\)), then $|K_n^{(q)}(x,y)|\le q^{a_n}$.
Moreover, if \(x,y\) lie in the same component of \(F\), then \(d_F(x,y)\ge d_n(x,y)>a_n\); otherwise \(L^{(q)}(x,y)=0\). Thus, $|L^{(q)}(x,y)|\le q^{a_n}$.

Combining the two cases, we obtain $|K_n^{(q)}(x,y)-L^{(q)}(x,y)|\le2q^{a_n}$ for any $x,y\in Z$. Since \(|Z|\le m\) and both matrices are self-adjoint, the Schur test yields
\[
\bigl\|K_n^{(q)}|_{Z\times Z}-L^{(q)}\bigr\| \le2mq^{a_n}.
\]
As \(L^{(q)}\ge0\), we conclude $K_n^{(q)}|_{Z\times Z}\ge-2mq^{a_n}\operatorname{Id}_{\ell^2(Z)}$. 
\end{proof}

\begin{proof}[Proof of Proposition~\ref{prop:maximal-norm-lower}]
It is clear from the assumption that $A=A^*\in\CC_u[X]$. Set $M:=\|A\|_{u,\max}$ and fix $\eta>0$.  Since $M \cdot \Id_{\ell^2(X)}-A$ is positive in $C^*_{u,\max}(X)$, \cite[Proposition~2.8]{MOSS} gives $B_1,\ldots,B_s\in\CC_u[X]$,
for some integer $s\ge1$, such that
\begin{equation}\label{eq:maximal-squares}
 (M+\eta)\cdot \Id_{\ell^2(X)}-A=\sum_{i=1}^s B_i^*B_i.
\end{equation}
Choose $r\ge0$ such that $\ppg(B_i)\le r$ for
$i=1,\ldots,s$, and set $m:=N(r)$.

Let $\{a_n\}_n$ be the sequence from Lemma~\ref{lem:signed-kernel} and for $q\in(0,1)$, set the kernel $K_n^{(q)}$ as above. Also set $\delta_n:=2mq^{a_n}$ and $J_n:=K_n^{(q)}+\delta_n \I_{\ell^2(Y_n)}$. By Lemma~\ref{lem:signed-kernel}, for all sufficiently large $n$ and $Z\subseteq Y_n$ with $|Z| \leq m$, we have $J_n|_{Z\times Z}\ge0$. Define a linear functional on $\B(\ell^2(X))$ by
\[
 \varphi_n(T)=\frac1{|Y_n|}\Tr(J_n\chi_{Y_n}T\chi_{Y_n}),
\]
where the trace is taken on $\ell^2(Y_n)$.

Claim: for each $i\in\{1,\ldots,s\}$, $\varphi_n(B_i^*B_i)\ge0$ for all sufficiently large $n$. Indeed, fix $i$ and simply denote $B:=B_i$. For each $z\in X$,  define 
\[
v_z:=\sum_{x\in Y_n}\overline{B_{z,x}}\delta_x
\in\ell^2(Y_n) \quad \text{where} \quad B_{z,x}:=\langle B\delta_x,\delta_z\rangle.
\]
Since $\ppg(B)\le r$, we have $\supp(v_z)\subseteq Y_n\cap B(z,r)$ and hence, $|\supp(v_z)|\le m$. Thus, $\langle J_nv_z,v_z\rangle\ge0$ for all sufficiently large $n$. Since
\[
(B^*B)_{x,y}=\sum_{z\in X}\overline{B_{z,x}}B_{z,y} \quad \text{for} \quad x,y\in Y_n,
\]
and $J_n$ is real symmetric, we obtain
\begin{align*}
|Y_n|\varphi_n(B^*B)&=\Tr(J_n\chi_{Y_n}B^*B\chi_{Y_n})=\sum_{x,y\in Y_n}J_n(y,x)(B^*B)_{x,y}=\sum_{z\in X}\sum_{x,y\in Y_n}
J_n(y,x)\overline{B_{z,x}}B_{z,y}\\
&=\sum_{z\in X}\sum_{x,y\in Y_n}
J_n(x,y)B_{z,y}\overline{B_{z,x}}=\sum_{z\in X}\langle J_nv_z,v_z\rangle\ge0.
\end{align*}
This proves the claim.

Since $K_n^{(q)}(x,x)=1$ for each $x\in Y_n$, $\varphi_n(\Id_{\ell^2(X)})=1+\delta_n$. Also $\chi_{Y_n}A\chi_{Y_n}=A_n$, which has zero diagonal.
For each edge $\{x,y\}\in E_n$, we have
$K_n^{(q)}(x,y)=q\varepsilon_n(\{x,y\})$ and
$(A_n)_{y,x}=\varepsilon_n(\{x,y\})$.
Hence, 
\begin{equation}\label{eq:maximal-A-value}
 \varphi_n(A)=\frac1{|Y_n|}\Tr(K_n^{(q)}A_n)
 =\frac{2q|E_n|}{|Y_n|}.
\end{equation}
Applying $\varphi_n$ to \eqref{eq:maximal-squares} and using \eqref{eq:maximal-A-value}, we obtain that for all sufficiently large $n$,
\[
 (M+\eta)(1+\delta_n)\ge
 \frac{2q|E_n|}{|Y_n|}.
\]
Since $m$ is fixed and $a_n\to\infty$, we have
$\delta_n=2mq^{a_n}\to0$. Then taking $\limsup$ as $n \to \infty$, letting $\eta\to 0$ and $q\to 1$, we conclude the proof.
\end{proof}

\subsection{Proof of Theorem~\ref{introthm:maximal-converse}}

We now combine the large-girth graph construction, the
signed-adjacency estimate and the maximal-norm lower bound to finish the proof of Theorem~\ref{introthm:maximal-converse}.

\begin{proof}[Proof of Theorem~\ref{introthm:maximal-converse}]
Suppose that $X$ fails Property~A.  Choose an integer $b\ge2$
sufficiently large such that $32\sqrt{2b}<b-2$.  We choose pairwise disjoint non-empty finite subsets $Y_n\subseteq X$ and the finite simple graphs $G_n=(Y_n, E_n)$ from Proposition~\ref{prop:large-girth-graphs}, which satisfy the conditions therein. Then Lemma~\ref{lem:signed-adjacency} provides signed adjacency
matrices $A_n$ on $Y_n$ such that $A:=\mathrm{(SOT)}$-$\sum_{n\in \NN}A_n$
belongs to $\CC_u[X]$ and satisfies
\[
 \|A\|=\sup_n\|A_n\|\le32\sqrt{2b}.
\]
On the other hand, Proposition~\ref{prop:maximal-norm-lower} gives
\[
 \|A\|_{u,\max}\ge
 \limsup_{n\to\infty}\frac{2|E_n|}{|Y_n|}\ge b-2.
\]
Thus $\|A\|<\|A\|_{u,\max}$, so $\lambda_X$ is not injective. By
Corollary~\ref{cor:scalar-reductions}(2),
$\lambda_{X,\H_0}$ is not injective for any nonzero separable
$\H_0$.

Conversely, if $X$ has Property~A, the result follows from \cite[Proposition~1.3]{SW}. Hence, we finish the proof.
\end{proof}

\subsection{Applications to discrete groups}

We conclude this section with applications of
Theorem~\ref{introthm:maximal-converse} to exactness
and crossed products of discrete groups.

Let $G$ be a countable discrete group. Denote by
$\lambda_g$ the left regular representation of
$g\in G$ on $\ell^2(G)$, and $u_g$ the canonical
unitary in $C^*_{\max}(G)$. Consider the $*$-algebra
\[
\mathcal A_G:=
\operatorname{span}\{
f\lambda_g\otimes u_g:
f\in\ell^\infty(G),\ g\in G
\}
\subseteq
\B(\ell^2(G))\odot C^*_{\max}(G),
\]
where $\B(\ell^2(G))\odot C^*_{\max}(G)$ is the algebraic tensor product. 
We can equip $\mathcal A_G$ with the restrictions of
the maximal and minimal tensor-product norms.

\begin{prop}\label{prop:LLP}
A countable discrete group $G$ is exact \emph{if and only if} the maximal
and minimal tensor-product norms agree on
$\mathcal A_G$.
\end{prop}

\begin{proof}
For a finite subset $F\subseteq G$ and
$f_g\in\ell^\infty(G)$ for $g\in F$, let
\[
a:=\sum_{g\in F}f_gv_g
\in\ell^\infty(G)\rtimes_{\mathrm{alg}}G
\quad \text{and} \quad 
z:=\sum_{g\in F}f_g\lambda_g\otimes u_g
\in\mathcal A_G,
\]
where $v_g$ denotes the canonical unitary in
the algebraic crossed product.

Let
$\iota:\ell^\infty(G)\hookrightarrow \B(\ell^2(G))$
be the equivariant inclusion, where $G$ acts
on $\B(\ell^2(G))$ by $\operatorname{Ad}\lambda$.
By \cite[Definition~3.1]{BEW}, for
$a\in\ell^\infty(G)\rtimes_{\mathrm{alg}}G$ we have
\[
\|a\|_{\mathrm{inj}}
\le
\|(\iota\rtimes G)(a)\|
\le
\|a\|_{\max}.
\]
Since \cite[Corollary~4.6]{BEW} gives
$\ell^\infty(G)\rtimes_{\mathrm{inj}}G
=\ell^\infty(G)\rtimes_{\max}G$,
both inequalities are equalities.
On the other hand, the standard untwisting
isomorphism for the inner action
$\operatorname{Ad}\lambda$
(see \cite[proof of Proposition~4.4]{BEW})
is given by
\[
\Phi:
\B(\ell^2(G))
\rtimes_{\operatorname{Ad}\lambda,\max}G
\longrightarrow
\B(\ell^2(G))
\otimes_{\max}C^*_{\max}(G),
\qquad
\Phi(Tv_g)=T\lambda_g\otimes u_g.
\]
In particular, we have $\Phi((\iota\rtimes G)(a))=z$. 
Consequently, we obtain
\[
\|z\|_{\max}
=\|(\iota\rtimes G)(a)\|
=\|a\|_{\ell^\infty(G)\rtimes_{\max}G}.
\]
On the other hand, Fell absorption gives
\[
\|z\|_{\min}
=
\left\|\sum_{g\in F}f_g\lambda_g
\right\|
=
\|a\|_{\ell^\infty(G)\rtimes_rG}.
\]

Equip $G$ with a proper right-invariant metric.
Then we have
\[
\ell^\infty(G)\rtimes_{\max}G
\cong C^*_{u,\max}(G)
\quad \text{and} \quad 
\ell^\infty(G)\rtimes_rG
\cong C^*_u(G).
\]
Consequently, the two tensor-product norms
agree on $\mathcal A_G$ if and only if
$C^*_{u,\max}(G)=C^*_u(G)$.
By Theorem~\ref{introthm:maximal-converse}
and the equivalence between Property~A and
exactness for discrete groups (see, \emph{e.g.}, \cite[Chapter 5]{BO}), we complete the proof.
\end{proof}

Finally, we prove Corollary~\ref{introcor:LLP}.

\begin{proof}[Proof of Corollary~\ref{introcor:LLP}]
We first assume that $G$ is countable. 

(1). Suppose that $C^*_{\max}(G)$ has the LLP. By Kirchberg's tensor-product characterisation (see \cite[Proposition~1.1]{Kirchberg1993} or \cite[Definition~4.3]{BEW}), the maximal and minimal tensor-product norms agree on $\B(\ell^2(G))\odot C^*_{\max}(G)$. In particular, they agree on $\mathcal A_G$, so Proposition~\ref{prop:LLP} implies that
$G$ is exact.

(2). Suppose that $\rtimes_{\mathrm{inj}}=\rtimes_r$. By \cite[Corollary~4.6]{BEW}, we have
\[
\ell^\infty(G)\rtimes_{\max}G
=
\ell^\infty(G)\rtimes_{\mathrm{inj}}G
=
\ell^\infty(G)\rtimes_rG.
\]
Since the maximal and reduced crossed products
of $\ell^\infty(G)$ are identified with
$C^*_{u,\max}(G)$ and $C^*_u(G)$, respectively,
Theorem~\ref{introthm:maximal-converse}
implies that $G$ has Property~A, and hence
is exact. The converse follows from
\cite[Proposition~4.2]{BEW}.

\medskip

Now let $G$ be an arbitrary discrete group.

(1). Suppose that $C^*_{\max}(G)$ has the LLP. By \cite[Remark~4.9]{BEW},
$C^*_{\max}(H)$ has the LLP for every countable
subgroup $H\leq G$. The countable case shows
that every such $H$ is exact. Since exactness
of discrete groups is determined by their
countable subgroups, it follows that $G$ is exact.

(2). Suppose that
$\rtimes_{\mathrm{inj},G}=\rtimes_{r,G}$.
We can deduce from \cite[Proposition~6.1]{BEW} that $\rtimes_{\mathrm{inj},H}
=
\rtimes_{r,H}$
for every countable subgroup $H\leq G$.
The countable case implies that every such
$H$ is exact, and hence $G$ is exact.
Conversely, if $G$ is exact, then
$\rtimes_{\mathrm{inj},G}=\rtimes_{r,G}$
by \cite[Proposition~4.2]{BEW}.
\end{proof}

\section{Strong quasi-locality}\label{sec:strong quasiloc}

This section is devoted to proving Theorem~\ref{introthm:sq=Roe}. Unlike Theorems~\ref{introthm:uniform-converse} and \ref{introthm:maximal-converse}, Theorem~\ref{introthm:sq=Roe} does \emph{not} involve Property~A. Our key ingredient is a Lipschitz transport lemma, which leads to a finite-propagation approximation estimate.

Note that Theorem~\ref{introthm:sq=Roe} holds trivially when $X=\emptyset$ or $\H_0=\{0\}$. Henceforth, we assume
$X\neq\emptyset$ and $\H_0\neq\{0\}$. The following is the Lipschitz transport result. 

\begin{lem}\label{lem:gram-transport}
Let $(X,d)$ be a discrete metric space of bounded geometry. Let $F$ be a finite subset of $X$, $E$ be a finite-dimensional Hilbert space, $S>0$ and $m\in \NN$. Suppose $a_{y,i}, b_{y,i}$ are vectors in $E$, where $y\in F$ and $1 \leq i \leq m$, satisfying:
\begin{equation}\label{eq:gram-agreement}
 \langle a_{y,i},a_{z,j}\rangle
 =\langle b_{y,i},b_{z,j}\rangle
 \quad\text{whenever} \quad y,z\in F \text{ with }d(y,z)\leq S \text{ and }\ 1\leq i,j\leq m.
\end{equation}
Then for each $x\in X$, there exists a contraction $\Psi_x\in\B(E)$ such that
\begin{equation*}
 \Psi_xa_{x,i}=b_{x,i}\quad \text{and} \quad \Psi_x^*b_{x,i}=a_{x,i} \quad \text{for} \quad x\in F \text{ and } 1 \leq i \leq m,
\end{equation*}
and the map $x\mapsto \Psi_x$ is $\frac{16}{S}$-Lipschitz, \emph{i.e.}, $\|\Psi_x-\Psi_z\|\leq\frac{16}{S}d(x,z)$ for any $x,z\in X$.
\end{lem}

\begin{proof}
For $x\in X$ and $0\leq t\leq S/2$, set
\[
V_x(t):=\operatorname{span}
\big\{a_{y,i}:y\in F \text{ with }\ d(x,y)\leq t \text{ and }\ 1\leq i\leq m\big\},
\]
and denote by $Q_x(t)$ the orthogonal projection
of $E$ onto $V_x(t)$. Clearly, we have $Q_x(t)\leq Q_x(t')$ whenever $0\leq t\leq t'\leq S/2$.

Fix $x\in X$. Since any two points in $B(x,S/2)$ have distance at most $S$, the assumption \eqref{eq:gram-agreement} implies that there is a unique isometry
\[
J_x:V_x(S/2)\longrightarrow E
\]
satisfying $J_xa_{y,i}=b_{y,i}$ for any $y\in F$ with $d(x,y)\leq S/2$ and $1\leq i\leq m$. 
Put $r:=S/8$ and define
\[
\Psi_x:=\frac{1}{r} \,J_x \int_r^{2r}Q_x(t)\,\mathrm{d}t.
\]
Since $F$ is finite, the integrand is a step function with finitely many discontinuities. Hence, the integral is well-defined in the
operator norm and clearly $\|\Psi_x\|\leq1$.

For $x\in F$ and $1\leq i\leq m$, we have
$Q_x(t)a_{x,i}=a_{x,i}$ for every
$t\in[r,2r]$. Consequently,
\[
\Psi_xa_{x,i}=J_xa_{x,i}=b_{x,i}.
\]
Since $J_x$ is an isometry and $J_xa_{x,i}=b_{x,i}$, we have $J_x^*b_{x,i}=a_{x,i}$. Hence,
\[
\Psi_x^*b_{x,i} =\frac1r\,\int_r^{2r}Q_x(t)J_x^*b_{x,i}\,\mathrm{d}t=\frac1r\,\int_r^{2r}Q_x(t)a_{x,i}\,\mathrm{d}t=a_{x,i}.
\]

It remains to prove the Lipschitz estimate. Fix $x,z\in X$ and put $h:=d(x,z)$. The case of $h=0$ is immediate. If $h\geq r$,
then
\begin{equation*}
\|\Psi_x-\Psi_z\|
\leq\|\Psi_x\|+\|\Psi_z\|
\leq2\leq\frac{16h}{S}.
\end{equation*}
Hence, we may assume $0<h<r$.
For any $t\in[r,2r]$, the triangle inequality gives
\begin{equation}\label{EQ:inclusion}
V_x(t-h)\subseteq V_z(t)\subseteq V_x(t+h)
\subseteq V_x(S/2).
\end{equation}
Since $J_x$ and $J_z$ agree on $V_z(t)$, we have $J_xQ_z(t)=J_zQ_z(t)$. It follows that
\[
\Psi_x-\Psi_z
=\frac1r J_x\int_r^{2r}
\bigl(Q_x(t)-Q_z(t)\bigr)\,dt.
\]

On the other hand, the inclusions in \eqref{EQ:inclusion} imply $Q_x(t-h)\leq Q_z(t)\leq Q_x(t+h)$. Hence for each $t\in [r,2r]$, we have
\begin{equation*}
Q_x(t)-Q_x(t+h) \le Q_x(t)-Q_z(t) \le Q_x(t)-Q_x(t-h).
\end{equation*}
Since $t\mapsto Q_x(t)$ is increasing, we have
\[
0 \leq\int_r^{2r} \bigl(Q_x(t+h)-Q_x(t)\bigr)\,\mathrm{d}t =\int_{2r}^{2r+h}Q_x(u)\,\mathrm{d}u -\int_r^{r+h}Q_x(u)\,\mathrm{d}u\leq h \cdot \Id_E.
\]
Similarly, we have
\[
0\leq\int_r^{2r} \bigl(Q_x(t)-Q_x(t-h)\bigr)\,\mathrm{d}t=\int_{2r-h}^{2r}Q_x(u)\,\mathrm{d}u-\int_{r-h}^{r}Q_x(u)\,\mathrm{d}u\leq h \cdot \Id_E.
\]
Combining the three inequalities above, we obtain
\[
-h \cdot \Id_E \leq\int_r^{2r} \bigl(Q_x(t)-Q_z(t)\bigr)\,\mathrm{d}t \leq h \cdot \Id_E.
\]
Therefore,
\[
\|\Psi_x-\Psi_z\| \leq\frac1r \left\|\int_r^{2r} \bigl(Q_x(t)-Q_z(t)\bigr)\,dt\right\| \leq\frac hr
=\frac{8h}{S} \leq\frac{16h}{S}.
\]
Together with the case $h\geq r$, this
proves the required estimate.
\end{proof}

Now we deal with the hard part of Theorem~\ref{introthm:sq=Roe}. We prove the following quantitative version.

\begin{prop}\label{prop:quantitative version}
Let $X$ be a discrete metric space of bounded geometry, and $\H_0$ be a separable Hilbert space. Given $\varepsilon, \delta, R>0$, suppose $T \in \B(\H_X)$ is self-adjoint and satisfies the following: for any map $g:X \rightarrow \KK(\K)_1$ with $(\delta, R)$-variation, $\big\| [T\otimes \Id_{\K} , \Lambda(g) ] \big\|< \varepsilon$. Then for any \(S>\frac{16R}{\delta}\), there exists a self-adjoint operator \(A\in\B(\H_X)\) with $\ppg(A) \leq S$ such that $\|T-A\| \leq \frac{\varepsilon}{2}$.
\end{prop}

\begin{proof}
For each $S>0$, denote $\MM_S:=\{A\in\B(\H_X): \ppg(A)\leq S\}$. We need to show that $\dist(T,\MM_S) \leq  \frac{\varepsilon}{2}$ for $S>\frac{16R}{\delta}$, where ``$\dist$'' denotes the operator norm distance.

Denote by $\S_1(\H_X)$ the set of trace-class operators on $\H_X$. For $D \in \S_1(\H_X)$, denote by $\Tr(D)$ the canonical trace and $\|D\|_1:=\Tr(|D|)$ the trace norm. Recall that $\B(\H_X)$ is the dual space of $\S_1(\H_X)$, given by $D \to \Tr(TD)$ for $D \in \S_1(\H_X)$ and $T \in \B(\H_X)$. Each $\MM_S$ is closed in the $w^*$-topology (also known as the ultraweak topology). Hence, the Hahn-Banach theorem for the $w^*$-topology on $\B(\H_X)$ shows:
\[
\dist(T,\MM_S) = \sup\left\{|\Tr(TD)|:
\begin{array}{l}
D\in \S_1(\H_X) \text{ with } \|D\|_1\leq 1 \\
\text{ and } \Tr(AD)=0 \text{ for all } A\in\MM_S
\end{array}
\right\}.
\]
For $D\in \S_1(\H_X)$, direct calculations show that $\Tr(AD)=0$ for all $A\in\MM_S$ if and only if $D_{x,y}=0$ whenever $d(x,y)\leq S$. Moreover, for any $D\in \S_1(\H_X)$ with $\|D\|_1\leq 1$ and $D_{x,y}=0$ whenever $d(x,y)\leq S$, choose $\lambda\in\CC$ with $|\lambda|=1$ such that $\lambda\Tr(TD)=|\Tr(TD)|$, and set $D'=\frac{\lambda D+\overline{\lambda}D^*}{2}$. Then $D'$ is self-adjoint, $\|D'\|_1\leq\|D\|_1$, and $D'_{x,y}=0$ whenever $d(x,y)\leq S$. Since $T=T^*$, we also have $\Tr(TD')=|\Tr(TD)|$. Consequently, we obtain
\begin{equation}\label{EQ:distance}
\dist(T,\MM_S) = \sup\left\{|\Tr(TD)|:
\begin{array}{l}
\text{self-adjoint } D\in\mathcal S_1(\H_X) \text{ with } \|D\|_1\leq 1\\
\text{ and } D_{x,y}=0 \text{ whenever } d(x,y)\leq S
\end{array}
\right\}.
\end{equation}

Let us fix $D\in\mathcal S_1(\H_X)$ satisfying the conditions in the right-hand side of \eqref{EQ:distance}. Choose increasing finite subsets $F_n\subseteq X$ with $\bigcup_{n\in \mathbb{N}}F_n=X$, and increasing finite-rank orthogonal projections $p_n\in\B(\H_0)$ converging strongly to $\Id_{\H_0}$. We may suppose $p_n\neq0$. Define
\[
 \Theta_n:=\chi_{F_n}(\Id_{\ell^2(X)}\otimes p_n)
 \quad \text{and} \quad D_n:=\Theta_nD\Theta_n.
\]
Since \(D \in \S_1(\H_X)\) and $\Theta_n\to\Id_{\H_X}$ in the strong operator topology, it follows that \(\|D_n-D\|_1\to 0\). Also $D_n$ is self-adjoint, $\|D_n\|_1\leq\|D\|_1 \leq 1$ and clearly, $(D_n)_{x,y}=0$ whenever $d(x,y)\leq S$. It therefore suffices to prove
\begin{equation}\label{EQ:Tr TDn}
|\Tr(TD_n)| \leq \frac{\varepsilon}{2} \|D_n\|_1 \quad \text{for each} \quad n\in \NN,
\end{equation}
since it implies $|\Tr(TD)|\leq \frac{\varepsilon}{2}$ and hence $\dist(T,\MM_S) \leq \frac{\varepsilon}{2}$.

Fix $n\in \NN$, and put
\[
 W_n:=\Theta_n\H_X=\ell^2(F_n)\otimes p_n\H_0
 \quad \text{and} \quad r_n:=\rank(p_n).
\]
Choose an orthonormal basis $\xi_1,\ldots,\xi_{r_n}$ for $p_n\H_0$ and put
\[
 e_{x,i}:=\delta_x\otimes\xi_i
 \quad \text{for} \quad x\in F_n \quad \text{and} \quad 1\leq i\leq r_n,
\]
where $\delta_x$ denotes the Dirac function at $x$. These vectors form an orthonormal basis of $W_n$.  
Regard $D_n$ as an operator on $W_n$ and decompose $D_n=D_n^+-D_n^-$ into its positive and negative parts. Since $(D_n)_{x,x}=0$ for each $x\in X$, we have $\Tr(D_n)=0$ and hence
\begin{equation}\label{eq:positive-traces}
 \Tr(D_n^+)=\Tr(D_n^-)=\frac12\|D_n\|_1.
\end{equation}

We aim to represent $\Tr(TD_n)$ as a matrix coefficient of the commutator in \eqref{EQ:commut strong quasi-local} for some $g$. Firstly, the positivity of $D_n^+$ and $D_n^-$ shows that the finite matrices
\[
 \bigl(\langle D_n^+ e_{y,j},e_{x,i}\rangle\bigr)_{(x,i),(y,j)} \quad \text{and} \quad \bigl(\langle D_n^- e_{y,j},e_{x,i}\rangle\bigr)_{(x,i),(y,j)}
\]
are positive semidefinite. Using Gram representations, there are vectors $a_{x,i},b_{x,i}\in W_n$ for $x\in F_n$ and $1\leq i\leq r_n$, satisfying
\begin{equation}\label{eq:positive-gram}
 \langle a_{x,i},a_{y,j}\rangle =\langle D_n^+e_{y,j},e_{x,i}\rangle \quad \text{and} \quad 
 \langle b_{x,i},b_{y,j}\rangle=\langle D_n^-e_{y,j},e_{x,i}\rangle
\end{equation}
for any $x,y\in F_n$ and $1\leq i,j\leq r_n$. 

Since $\K$ is infinite-dimensional, choose an isometric embedding $\iota_n:W_n\to\K$. Now we set
\begin{equation*}
 u_n:=\sum_{x\in F_n}\sum_{i=1}^{r_n}e_{x,i}\otimes\iota_na_{x,i} 
 \quad \text{and} \quad
 v_n:=\sum_{x\in F_n}\sum_{i=1}^{r_n}e_{x,i}\otimes\iota_nb_{x,i}.
\end{equation*}
By \eqref{eq:positive-traces} and \eqref{eq:positive-gram}, we have
\begin{equation}\label{eq:uv-norms}
 \|u_n\|^2=\Tr(D_n^+)=\frac12\|D_n\|_1
 \quad \text{and} \quad \|v_n\|^2=\Tr(D_n^-)=\frac12\|D_n\|_1.
\end{equation}
Moreover, using \eqref{eq:positive-gram} we obtain
\begin{align}\label{eq:trace-form-plus}
\langle(T\otimes\Id_{\K})u_n,u_n\rangle & = \sum_{x,y\in F_n}\sum_{i,j=1}^{r_n} \langle Te_{x,i}\otimes\iota_na_{x,i}, e_{y,j} \otimes \iota_na_{y,j} \rangle \\
&=\sum_{x,y\in F_n}\sum_{i,j=1}^{r_n} \langle Te_{x,i},e_{y,j}\rangle \langle D_n^+e_{y,j},e_{x,i}\rangle =\Tr(TD_n^+). \notag
\end{align}
Similarly, we have
\begin{equation}\label{eq:trace-form-minus}
\langle(T\otimes\Id_{\K})v_n,v_n\rangle = \Tr(TD_n^-).
\end{equation}

Since $(D_n)_{x,y}=0$ when $d(x,y)\leq S$, it follows from \eqref{eq:positive-gram} that
\[
 \langle a_{x,i},a_{y,j}\rangle
 =\langle b_{x,i},b_{y,j}\rangle
 \quad\text{for} \quad x,y\in F_n \text{ with } d(x,y)\leq S \text{ and } 1\leq i,j\leq r_n.
\]
Applying Lemma~\ref{lem:gram-transport} to $E=W_n$, $F=F_n$ and $m=r_n$, we obtain contractions $\Psi_x\in\B(W_n)$ for each $x\in X$, such that
\begin{equation}\label{EQ:transport identity}
 \Psi_xa_{x,i}=b_{x,i}\quad \text{and} \quad \Psi_x^*b_{x,i}=a_{x,i} \quad \text{for} \quad x\in F_n \text{ and } 1\leq i\leq r_n,
\end{equation}
and 
\[
 \|\Psi_x-\Psi_y\|\leq\frac{16}{S}d(x,y) \quad \text{for any} \quad x,y\in X.
\]
Define $g_n: X \to \KK(\K)$ by $g_n(x):=\iota_n\Psi_x\iota_n^*$ for $x\in X$. Then each $g_n(x)$ has finite rank and norm at most one. Moreover, for $x,y\in X$ with $d(x,y)< R$, we have
\[
 \|g_n(x)-g_n(y)\|\leq\frac{16}{S}d(x,y)<\delta.
\]
Thus, $g_n$ has $(\delta, R)$-variation. Moreover, \eqref{EQ:transport identity} shows:
\begin{equation}\label{eq:uv-transport}
 \Lambda(g_n)u_n=v_n \quad \text{and} \quad \Lambda(g_n)^*v_n=u_n.
\end{equation}
Combining \eqref{eq:trace-form-plus}, \eqref{eq:trace-form-minus} and \eqref{eq:uv-transport}, we obtain
\begin{align*}
 \langle[T\otimes\Id_{\K},\Lambda(g_n)]u_n,v_n\rangle
 &=\langle(T\otimes\Id_{\K})v_n,v_n\rangle
   -\langle(T\otimes\Id_{\K})u_n,u_n\rangle\\
 &=\Tr(TD_n^-)-\Tr(TD_n^+)=-\Tr(TD_n).
\end{align*}
Hence, the assumption together with \eqref{eq:uv-norms} gives
\[
 |\Tr(TD_n)|\leq \big\|[T\otimes\Id_{\K},\Lambda(g_n)]\big\| \cdot \|u_n\| \cdot \|v_n\|
 \leq\frac{\varepsilon}{2} \cdot \|D_n\|_1.
\]
This proves \eqref{EQ:Tr TDn}, and \eqref{EQ:distance} gives $\dist(T,\MM_S) \leq \frac{\varepsilon}{2}$.

Since $\MM_S$ is weak-\(\ast\) closed, its intersection with any norm-closed ball is weak-\(\ast\) compact by the Banach–Alaoglu theorem. Hence there exists \(A\in\MM_S\) such that $\|T-A\|=\dist(T,\MM_S) \le\frac{\varepsilon}{2}$. Since $T$ is self-adjoint, replacing $A$ by $\frac{A+A^*}{2}$ gives a self-adjoint operator satisfying the same estimate. This completes the proof.
\end{proof}

Now we are ready to prove Theorem~\ref{introthm:sq=Roe}.

\begin{proof}[Proof of Theorem~\ref{introthm:sq=Roe}]
Firstly, suppose that $T\in\B(\H_X)$ is strongly quasi-local. Since
\[
[T^*\otimes\Id_{\K},\Lambda(g)]
=-[T\otimes\Id_{\K},\Lambda(g^*)]^*,
\]
the operators $T_1:=\frac{T+T^*}{2}$ and $T_2:=\frac{T-T^*}{2i}$ are strongly quasi-local and self-adjoint.

For $\varepsilon>0$, choose $\delta,R>0$ such that both $T_1$ and $T_2$ satisfy the assumption of Proposition~\ref{prop:quantitative version} with these parameters. Choose $S>\frac{16R}{\delta}$. By Proposition~\ref{prop:quantitative version}, there exist self-adjoint operators $A_1,A_2\in\B(\H_X)$ with propagation at most $S$ such that
\[
\|T_1-A_1\|\leq\frac{\varepsilon}{2}
\quad \text{and} \quad
\|T_2-A_2\|\leq\frac{\varepsilon}{2}.
\]
Setting $A:=A_1+iA_2$, we obtain $\ppg(A)\leq S$ and 
\[
\|T-A\| \leq\|T_1-A_1\|+\|T_2-A_2\| \leq\varepsilon.
\]

Conversely, the argument in the proof of \cite[Proposition~3.6(2)]{BCZ23} shows that every bounded finite-propagation operator is strongly quasi-local. Since strong quasi-locality is preserved under norm limits, the converse follows.

Finally, we prove the ``Furthermore'' part. Fix $T \in C^*_{sq}(X;\H_0)$ and $\varepsilon>0$. By the first part, we can choose a finite-propagation operator $B\in\B(\H_X)$ such that $\|T-B\|<\frac{\varepsilon}{2}$. For each $x\in X$, denote $\H_x:= \chi_{\{x\}} \H_X \cong \H_0$, and choose positive numbers $\varepsilon_x$ such that
\[
\sum_{x\in X}\varepsilon_x^2
<\frac{\varepsilon^2}{16}.
\]
Since $\chi_{\{x\}}T$ and $\chi_{\{x\}}T^*$ are compact for each $x\in X$, there exists a finite-rank orthogonal projection $P_x$ on $\H_x$ such that
\[
\|(\Id_{\H_x}-P_x)\chi_{\{x\}}T\|
\leq\varepsilon_x 
\quad \text{and} \quad 
\|(\Id_{\H_x}-P_x)\chi_{\{x\}}T^*\|
\leq\varepsilon_x.
\]
Set $P:=\mathrm{(SOT)}$-$\sum_{x\in X}P_x$. Then $\|(\Id_{\H_X}-P)T\|<\frac{\varepsilon}{4}$ and $\|T(\Id_{\H_X}-P)\|<\frac{\varepsilon}{4}$, which shows that
\[
\|T-PTP\|<\frac{\varepsilon}{2}.
\]

Set $A:=PBP$. Clearly, $A$ is locally compact and has finite propagation. Moreover, we have
\[
\|T-A\| \leq\|T-PTP\|+\|P(T-B)P\| \leq\|T-PTP\|+\|T-B\| <\varepsilon.
\]
Thus $T\in C^*(X;\H_0)$. The reverse inclusion is exactly \cite[Proposition~3.6(2)]{BCZ23}.
\end{proof}

\begin{rem}
The proof of Theorem~\ref{introthm:sq=Roe} shows that strong quasi-locality can be characterised using only finite-dimensional matrix-valued test functions $g: X \to M_m(\CC)_1$, provided that the parameters \(\delta,R\) are chosen independently of \(m\in \NN\). On the other hand, \cite[Theorem 2.8]{ST19} shows that ordinary quasi-locality can be characterised by scalar-valued test functions, or equivalently by matrix-valued test functions of each fixed size \(m\), where \(\delta,R\) may depend on \(m\). 
\end{rem}

\bibliographystyle{plain}
\bibliography{quasi_local_and_WCP}

\begin{thebibliography}{10}

\bibitem{AGS}
Goulnara Arzhantseva, Erik Guentner, and J{\'a}n {\v{S}}pakula.
\newblock Coarse non-amenability and coarse embeddings.
\newblock {\em Geom. Funct. Anal.}, 22(1):22--36, 2012.

\bibitem{BCZ23}
Hengda Bao, Xiaoman Chen, and Jiawen Zhang.
\newblock Strongly quasi-local algebras and their {$K$}-theories.
\newblock {\em J. Noncommut. Geom.}, 17(1):241--285, 2023.

\bibitem{BNSWW}
Jacek Brodzki, Graham~A. Niblo, J{\'a}n {\v S}pakula, Rufus Willett, and Nick
  Wright.
\newblock Uniform local amenability.
\newblock {\em J. Noncommut. Geom.}, 7(2):583--603, 2013.

\bibitem{BO}
Nathanial~P. Brown and Narutaka Ozawa.
\newblock {\em {$C^*$}-Algebras and Finite-Dimensional Approximations},
  volume~88 of {\em Graduate Studies in Mathematics}.
\newblock American Mathematical Society, 2008.

\bibitem{BEW}
Alcides Buss, Siegfried Echterhoff, and Rufus Willett.
\newblock The maximal injective crossed product.
\newblock {\em Ergodic Theory Dynam. Systems}, 40(11):2995--3014, 2020.

\bibitem{CGZ24}
Xiaoman Chen, Kun Gao, and Jiawen Zhang.
\newblock The strongly quasi-local coarse {N}ovikov conjecture and {B}anach
  spaces with {P}roperty {(H)}.
\newblock {\em Kyoto J. Math.}, 64(1):31--74, 2024.

\bibitem{Elek}
G{\'a}bor Elek.
\newblock Uniform local amenability implies property {A}.
\newblock {\em Proc. Amer. Math. Soc.}, 149(6):2573--2577, 2021.

\bibitem{ET}
G{\'a}bor Elek and {\'A}d{\'a}m Tim{\'a}r.
\newblock Strong almost finiteness.
\newblock {\em J. Funct. Anal.}, 289(10):111116, 2025.

\bibitem{Eng18}
Alexander Engel.
\newblock Index theorems for uniformly elliptic operators.
\newblock {\em New York J. Math.}, 24:543--587, 2018.

\bibitem{GWY08}
Guihua Gong, Qin Wang, and Guoliang Yu.
\newblock Geometrization of the strong {N}ovikov conjecture for residually
  finite groups.
\newblock {\em J. Reine Angew. Math.}, 621:159--189, 2008.

\bibitem{GMR18}
Rostislav Grigorchuk, Magdalena Musat, and Mikael R{\o}rdam.
\newblock Just-infinite {$C^*$}-algebras.
\newblock {\em Comment. Math. Helv.}, 93(1):157--201, 2018.

\bibitem{GK02}
Erik Guentner and Jerome Kaminker.
\newblock Exactness and the {N}ovikov conjecture.
\newblock {\em Topology}, 41(2):411--418, 2002.

\bibitem{JZZ23}
Baojie Jiang, Jiawen Zhang, and Jianguo Zhang.
\newblock Quasi-locality for {\'e}tale groupoids.
\newblock {\em Comm. Math. Phys.}, 403:329--379, 2023.

\bibitem{Kirchberg1993}
Eberhard Kirchberg.
\newblock On nonsemisplit extensions, tensor products and exactness of group
  {$C^*$}-algebras.
\newblock {\em Invent. Math.}, 112(3):449--489, 1993.

\bibitem{KLS}
Rasmus Kyng, Kyle Luh, and Zhao Song.
\newblock Four deviations suffice for rank 1 matrices.
\newblock {\em Adv. Math.}, 375:107366, 2020.

\bibitem{LNSZ21}
Kang Li, Piotr~W. Nowak, J{\'a}n {\v{S}}pakula, and Jiawen Zhang.
\newblock Quasi-local algebras and asymptotic expanders.
\newblock {\em Groups Geom. Dyn.}, 15(2):655--682, 2021.

\bibitem{LZZ26}
Kang Li, Jiawen Zhang, and Jingming Zhu.
\newblock {$K$}-theoretic comparison of {R}oe and quasi-local algebras via
  projections, 2026.
\newblock arXiv:2608.22439.

\bibitem{McD}
Colin McDiarmid.
\newblock On the method of bounded differences.
\newblock In Johannes Siemons, editor, {\em Surveys in Combinatorics, 1989},
  volume 141 of {\em London Mathematical Society Lecture Note Series}, pages
  148--188. Cambridge University Press, 1989.

\bibitem{MOSS}
Masato Mimura, Narutaka Ozawa, Hiroki Sako, and Yuhei Suzuki.
\newblock Group approximation in {Cayley} topology and coarse geometry, {III}:
  Geometric property ({T}).
\newblock {\em Algebr. Geom. Topol.}, 15(2):1067--1091, 2015.

\bibitem{Ozawa2000}
Narutaka Ozawa.
\newblock Amenable actions and exactness for discrete groups.
\newblock {\em C. R. Acad. Sci. Paris S\'er. I Math.}, 330(8):691--695, 2000.

\bibitem{Oza}
Narutaka Ozawa.
\newblock Embeddings of matrix algebras into uniform {Roe} algebras and
  quasi-local algebras.
\newblock {\em J. Eur. Math. Soc.}, 2025.
\newblock Published online first.

\bibitem{Roe88}
John Roe.
\newblock An index theorem on open manifolds. {I}, {II}.
\newblock {\em J. Differential Geom.}, 27(1):87--113, 115--136, 1988.

\bibitem{Roe93}
John Roe.
\newblock Coarse cohomology and index theory on complete {Riemannian}
  manifolds.
\newblock {\em Mem. Amer. Math. Soc.}, 104(497):x+90, 1993.

\bibitem{Roe96}
John Roe.
\newblock {\em Index theory, coarse geometry, and topology of manifolds},
  volume~90 of {\em CBMS Regional Conference Series in Mathematics}.
\newblock Published for the Conference Board of the Mathematical Sciences,
  Washington, DC; by the American Mathematical Society, Providence, RI, 1996.

\bibitem{RW}
John Roe and Rufus Willett.
\newblock Ghostbusting and property {A}.
\newblock {\em J. Funct. Anal.}, 266(3):1674--1684, 2014.

\bibitem{Sako}
Hiroki Sako.
\newblock A generalization of expander graphs and local reflexivity of uniform
  {Roe} algebras.
\newblock {\em J. Funct. Anal.}, 265(7):1367--1391, 2013.

\bibitem{Sch}
Isaac~J. Schoenberg.
\newblock Metric spaces and positive definite functions.
\newblock {\em Trans. Amer. Math. Soc.}, 44:522--536, 1938.

\bibitem{STY02}
Georges Skandalis, Jean-Louis Tu, and Guoliang Yu.
\newblock The coarse {Baum-Connes} conjecture and groupoids.
\newblock {\em Topology}, 41(4):807--834, 2002.

\bibitem{ST19}
J{\'a}n {\v{S}}pakula and Aaron Tikuisis.
\newblock Relative commutant pictures of {Roe} algebras.
\newblock {\em Comm. Math. Phys.}, 365(3):1019--1048, 2019.

\bibitem{SW}
J{\'a}n {\v S}pakula and Rufus Willett.
\newblock Maximal and reduced {Roe} algebras of coarsely embeddable spaces.
\newblock {\em J. Reine Angew. Math.}, 678:35--68, 2013.

\bibitem{SZ20}
J{\'a}n {\v{S}}pakula and Jiawen Zhang.
\newblock Quasi-locality and {P}roperty {A}.
\newblock {\em J. Funct. Anal.}, 278(1):108299, 2020.

\bibitem{Willett}
Rufus Willett.
\newblock A non-amenable groupoid whose maximal and reduced {$C^*$}-algebras
  are the same.
\newblock {\em M{\"u}nster J. Math.}, 8:241--252, 2015.

\bibitem{WillettYu2020}
Rufus Willett and Guoliang Yu.
\newblock {\em Higher index theory}, volume 189 of {\em Cambridge Studies in
  Advanced Mathematics}.
\newblock Cambridge University Press, Cambridge, 2020.

\bibitem{Yu00}
Guoliang Yu.
\newblock The coarse {B}aum-{C}onnes conjecture for spaces which admit a
  uniform embedding into {H}ilbert space.
\newblock {\em Invent. Math.}, 139(1):201--240, 2000.

\end{thebibliography}

\end{document}